\documentclass[12pt,a4paper]{article}
\usepackage[T1]{fontenc}
\usepackage{lmodern}
\usepackage[latin1]{inputenc}
\usepackage{latexsym}                 
\usepackage{color}                 
\usepackage{amssymb}
\usepackage{amsmath}
\usepackage{amsthm}
\usepackage{graphicx}
\usepackage{hyperref}
\usepackage{verbatim}
\usepackage{mathptmx}      % use Times fonts if available on your TeX system
\def\RSfootnote#1{%
}
\def\eref#1{(\ref{#1}%
)}
\def\RSref#1{\ref{#1}%
}
\def\RSlabel#1{\label{#1}%
}
\def\RScite#1{\cite{#1}%
}

\newcommand{\bql}[1]{%
\begin{equation}\label{#1}%
}

\def\filename#1{}

\newcommand{\eq}{\end{equation}}
\def\dfrac#1#2{\displaystyle{\frac{#1}{#2}   }}

\def\b1{\mathbf 1}

\def\biglf{\par\bigskip\noindent}
\def\biglf{\par\bigskip\noindent}
\begin{document}
\begin{center}
  {\bf A Simple Solution for Maximum Range Flight
    %Fuel-to-distance-optimal quasi-steady flight
  }
\biglf 
Robert Schaback\footnote{Prof. Dr. R. Schaback\\
  Institut für Numerische und
  Angewandte Mathematik, \\
  Lotzestraße 16-18, D-37083 Göttingen, Germany\\
  schaback@math.uni-goettingen.de\\
http://www.num.math.uni-goettingen.de/schaback/ } \\
\biglf
Draft of \today
\end{center}
    {\bf Abstract}: Within the standard framework of quasi-steady flight,
    this paper derives a speed that realizes the maximal
    obtainable range per unit of  fuel. If this speed is chosen
    at each instant of
    a flight plan $h(x)$ giving altitude $h$ as a function of distance $x$,
    a variational problem for finding an optimal $h(x)$ can be formulated
    and solved. It yields flight plans with maximal range, and
    these turn out to consist of mainly three phases using the optimal speed:
    starting with a climb at maximal
    continuous admissible
    thrust, ending with a continuous descent
    at idle thrust, and in between with a transition based
    on a solution of the Euler-Lagrange
    equation for the variational problem. A similar variational problem
    is derived and solved for speed-restricted flights,
    e.g. at 250 KIAS below 10000 ft. In contrast to the literature,
    the approach of this paper does not need more than
    standard ordinary differential equations solving variational problems
    to derive range-optimal trajectories.
Various numerical examples
    based on a Standard Business Jet are added for illustration.
    %****************************************************************
%%%%%%%%%%%%%%%%%%%%%%%%%%%%%%%%%%%%%%%%%%%%%%%%%%%%%%%%%
    \section{Introduction}\RSlabel{SecIntro}

% \red{\RScite{dalmau-prats:2014-1} global, optimal control, realistic}

    The problem of calculating
    flight trajectories that minimize fuel consumption
    or maximize range 
    has a long history,
    see e.g. the references in \RScite{peckham:1974-1,vinh:1980-1,%
      burrows:1982-1, pierson-ong:1989-1,
      valenzuela-rivas:2014-1}. 
    Various mathematical techniques were applied, ranging from
    energy considerations
    \RScite{rutowski:1954-1,bryson-et-al:1969-1,calise:1977-1,yajnik:1977-1}, 
    parametrizations of trajectories
    \RScite{rader-hull:1975-1, burrows:1982-1,valenzuela-rivas:2014-1} via
    certain forms of Optimal
    Control Theory
    \RScite{franco-rivas:2014-1, garciaheras-et-al:2016-1,park-clarke:2015-1}
    to Multiobjective Optimization using various cost functionals
    \RScite{maazoun:2015-1,gardi-et-al:2016-1, saucier-et-al:2017-1}.
    Compilations of numerical methods for trajectory
    calculation and optimization are in
    \RScite{betts:1998-1,huang-et-al:2012-1}.
    \biglf
    A particularly simple
    solution for range-optimal flight is well-known in case of
    {\em horizontal} flight, see e.g. \RScite{peckham:1974-1,vinh:1995-1,%
    stengel:2004-1,myose-et-al:2005-1}.
    It follows from maximizing
    the ratio $\sqrt{C_L}/C_D$ of the lift and drag coefficients,
    leading to a speed that is by a factor $\sqrt[4]{3}=1.316$ larger than
    the speed maximizing the lift-to-drag ratio ${C_L}/C_D$.
    This paper provides an extension to general non-horizontal flight,
    staying close to basic classroom texts
    \RScite{miele:1962-1,vinh:1980-1,vinh:1995-1,stengel:2004-1,%
      hull:2007-1,phillips:2010-1}  and focusing on
    standard numerical methods that just
    solve systems of ordinary differential equations. There is no constraint on
    fixed altitude, but wind effects and fixed arrival times are ignored
    \RScite{franco-et-al:2010-1,franco-rivas:2014-1}.
    \biglf
    Starting with the basics of {\em quasi-steady flight} in Section
    \RSref{SecQSF} and an arbitrary given {\em flight path}
    in terms of a function $h(x)$ of altitude $h$
    of distance $x$,
    a specific speed assignment that maximizes range at each
    instant of the flight is calculated in Section \RSref{SecRMaxHF}.
    Then Section \RSref{SecVP} varies
    flight paths with range-optimal speed assignments and derives
    a {\em variational} problem that gets range-optimal flight
    trajectories by solving a second-order
    Euler-Lagrange differential equation for
    $h(x)$. But the solutions may violate thrust restrictions. Therefore
    the variational problem is a {\em constrained} one, and its solutions
    must either satisfy the Euler-Lagrange equation or follow one of the
    restrictions. Section \RSref{SecCROT} provides the solutions for
    thrust-constrained maximal range trajectories, and these occur
    for climb/cruise at maximal continuous admissible thrust
    and for Continuous Descent at idle thrust, both still using the
    speed assignment of Section \RSref{SecRMaxHF}. Between these
    two range-optimal
    trajectory parts, the {\em transition} from maximal to idle
    thrust must follow the Euler-Lagrange equation of Section
    \RSref{SecVP}, solving the range optimality problem completely 
    for flights above 10000 ft.
    \biglf
    Below 10000 ft, the speed restriction
    to 250 knots indicated airspeed (KIAS) comes into play.
    Since the unconstrained solutions of the variational problem
    severely violate the speed restriction, a range-optimal
    solution has to follow the 250 KIAS restriction below 10000 ft.
    Therefore a second variational
    problem is derived in Section \RSref{SecPS} that allows to
    calculate range-optimal trajectories under speed restriction,
    and the outcome is similar to the previous situation. Optimal trajectories
    violate thrust restrictions, and thus they either follow a thrust
    restriction or satisfy a second Euler-Lagrange equation.
    The result is that a range-maximal climb strategy below 10000 ft
    at 250 KIAS 
    first uses maximal admissible thrust and then continues
    with a solution of the second Euler-Lagrange equation.
    Since all of this ignores restrictions by Air Traffic Control,
    Section \RSref{SecQSFFLC} deals with flight level changes
    between level flight sections at range-optimal speed. All
    trajectory parts derived so far are combined by the final Section
    \RSref{SecFPfOFU}.
    \biglf
    The mathematical procedures to calculate range-optimal 
    trajectories are simple enough to be carried out rather quickly
    by any reasonably fast 
    and suitably programmed Flight Management System, and
    the 
    range-optimal speed could be displayed on any Electronic
    Flight Instrument System.
    \biglf
    All model calculations were done for the Standard Business Jet (SBJ)
    of \RScite{hull:2007-1} for convenience, using
    the simple turbojet propulsion model presented there. 
    Symbolic formula manipulations,
    e.g. for setting up the Euler equations for the two variational problems,
    were done by MAPLE$^\copyright$, and MATLAB$^\copyright$
    was used for all numerical  calculations, mainly ODE solving. 
    Programs are available from the author on request. 
%\RScite{burrows:1982-1, franco-rivas:2015-1, garciaheras-et-al:2016-1,
%  park-clarke:2015-1,valenzuela-rivas:2014-1, saucier-et-al:2017-1,
%  gardi-et-al:2016-1, pierson-ong:1989-1,
 % % murrietamendoza-botez:2015-1,
%  franco-rivas:2014-1,franco-et-al:2010-1,maazoun:2015-1}
%%%%%%%%%%%%%%%%%%%%%%%%%%%%%%%%%%%%%%%%%%%%%%%%%%%%%%%%%
\section{Quasi-Steady Flight}\RSlabel{SecQSF}
The standard \RScite{vinh:1980-1,%
  vinh:1995-1,stengel:2004-1,hull:2007-1, phillips:2010-1}
equations of {\em quasi-steady flight} are 
\bql{eqbas2Dc}
\begin{array}{rcl}
  \dot{x}&=& V\cos \gamma\\
  \dot{h} &=& V\sin \gamma\\
  \dot{W} &=& -C\,T\\
  0&=&
  T-D-W\sin \gamma\\[0.2cm]
  0&=&L-W\cos \gamma.
\end{array}
\eq
with {\em distance} $x$, {\em altitude} $h$, {\em true airspeed} $V$,
{\em flight path angle} $\gamma$,
{\em specific fuel consumption} $C$, {\em weight} $W$, {\em thrust} $T$,
{\em drag} $D$, and {\em lift} $L$.
Like weight, lift, and drag, we consider
thrust as a force, not a mass. Furthermore, we omit the influence of flaps,
spoilers, or extended gears, i.e. we exclusively work in
{\em clean configuration}. The equations live on short time intervals
where speed $V$ and angle $\gamma$ are considered to be constant,
but they will lead to useful equations that 
describe long-term changes of $V$ and $\gamma$. Throughout the paper,
we shall assume that the specific fuel consumption $C$ is independent of speed,
but dependent on altitude $h$.
\biglf
    {\em Lift} and {\em Drag} are 
$$
L=\dfrac{1}{2}C_L\rho V^2 S,\;D=\dfrac{1}{2}C_D\rho V^2 S
$$
with the altitude-dependent {\em air density} $\rho$, the
{\em wing planform area} $S$, and the specific lift and drag coefficients
$C_L$ and $C_R$. We also use the 
{\em drag polar}
\bql{eqdragpolar}
C_D=C_{D_0}+KC_L^2
\eq
for further analysis.
The {\em induced drag factor} $K$
and the {\em lift-independent drag coefficient}   
$C_{D_0}$ are dependent on Mach number,
but we ignore this fact for simplicity.
When it comes to calculations, and if speed
and altitude are known, one can insert the Mach-dependent values whenever
necessary, but we did not implement this feature
and completely ignore the minor
dependence on Reynolds number and viscosity. Throughout, we shall use
the well-known exponential model
\bql{eqairdens}
\rho(h)=1.225\exp(-h/9042)
\eq
for air density in $kg/m^3$  as a function of altitude $h$ in $m$.
\biglf
If $W,\,\gamma,\,x$, and $h$
are considered to be independent
variables and $S,\,C_{D_0},\,K$ to be constants,
we have five equations for the six unknowns
$T,\,D,\,L,\,C_D,\,C_L,\,V$, leaving one variable for optimization
that we are free to choose. Whatever will be optimized later, the solution will
not depend on the choice of the remaining variable. Because pilots can
fly prescribed speeds or prescribed thrusts
and cannot directly maintain certain values of $C_L$,
there is a certain practical preference for $V$ and $T$.
Both of the latter are restricted in practice, and these restrictions
will need special treatment.
\biglf
Because many following calculations will be simpler,
we introduce
$$
R:=\dfrac{1}{2}\rho V^2 \dfrac{S}{W}=\dfrac{\dfrac{1}{2}\rho V^2}{\dfrac{W}{S}}
$$
as the ratio between {\em dynamic pressure} $\frac{1}{2}\rho V^2$ and
{\em wing pressure} $W/S$ and call it the {\em pressure ratio}.
Avoiding mass notions, we prefer {\em wing pressure} over the usual {\em wing
  loading}. It will turn out that the  pressure ratio $R$
is of central importance when dealing with quasi-steady flight.
It combines speed, altitude (via $\rho$), weight, and wing planform area
into a very useful dimensionless quantity that should get more
attention by standard texts on Flight Mechanics.
The variable $R$ arises in \cite[p. 201, (3)]{miele:1962-1}
temporarily, as some $u^2$, and later (p. 216)
as $M^2/\omega$ in various expressions
where $M$ is Mach number and $\omega$ is the dimensionless wing loading
$$
\omega=\dfrac{2W}{\rho S a^2}
$$
with $a$ being the speed of sound. This $M^2/\omega$
coincides with our $R$ written
in terms of Mach number instead of true airspeed. 
\biglf
Here, we use $R$ to express the other possible variables via $R$,
and then each pair of the variables can be connected via $R$.
The results are
\bql{eqRconv} 
\begin{array}{rcl}
  C_L&=& \dfrac{\cos \gamma}{R},\;
 C_D=C_{D_0}+\dfrac{K}{R^2}\cos^2\gamma,\\[0.3cm]
\dfrac{L}{W}&=&\cos\gamma,\;
\dfrac{D}{W}=C_{D_0}R+\dfrac{K}{R}\cos^2\gamma,\\[0.3cm]
   V^2&=&\dfrac{2RW}{\rho S}
\end{array}
\eq
and in particular
\bql{eqTWfull}
   \dfrac{T}{W}
  =
  C_{D_0}R+
  \dfrac{K\cos^2 \gamma}{R}+\sin \gamma
\eq
after some simple calculations.
By taking thrust, lift, and drag {\em relative} to the current weight,
their equations become dimensionless. The advantage is that the weight
drops out of many of the following arguments. Only the constants $K$ and
$C_{D_0}$ of the drag polar are relevant.
\biglf
Various texts, e.g. \RScite{phillips:2010-1},
reduce everything to $C_L$. Because of $C_L=1/R$ 
for horizontal flight, this is not much different from working with $R$,
but for general flight angles it will pay off to work with $R$.
\biglf
To illustrate how $R$ resembles speed but hides altitude and weight,
consider horizontal flight
and the maximal $C_L$ value $C_L^{stall}$ that belongs to the maximal angle
of attack. Then $C_L= \frac{1}{R}$ shows that $R_{stall}=1/C_L^{stall}$
is the minimal admissible $R$, usually around 0.8. Then the stall speed
as a function of weight, altitude, and wing loading  is
$V^2_{stall}=\dfrac{2R_{stall}W}{\rho S}$, derived from the constant $R_{stall}$.
\biglf
As another example,
consider the usual
maximization of the lift-to-drag
ratio $L/D$. This means minimization of the denominator of
$$
\dfrac{L}{D}=\dfrac{\cos\gamma}{RC_{D_0}+\dfrac{K}{R}\cos^2\gamma},
$$
with respect to $R$, leading to 
\bql{eqRLD}
R^2_{L/D}=\dfrac{K}{C_{D_0}}\cos^2\gamma.
\eq
From here, the other variables follow via \eref{eqRconv}, e.g.
$$
\begin{array}{rcl}
V^2_{L/D}&=&\dfrac{2R_{L/D}W}{\rho S}
=\dfrac{2W\sqrt{K}}{\rho S\sqrt{C_{D_0}}}\cos\gamma.
%\dfrac{T_{L/D}}{W}
%  &=&2\sqrt{KC_{D_0}}\cos\gamma+\sin\gamma.
\end{array} 
$$
From the equation
$$
\dfrac{T}{W}-\sin\gamma
=
C_{D_0}R+
\dfrac{K\cos^2 \gamma}{R}=\dfrac{D}{W}=\dfrac{D}{L}\cos\gamma
$$
it follows that this solution realizes the minimal $T/W$ ratio for given
$\gamma$ as well, i.e. it is the solution for minimal thrust.
This implies the inequality
\bql{eqTWrestr}
\dfrac{T}{W}-\sin\gamma\geq2\sqrt{KC_{D_0}}\cos\gamma=
\dfrac{T_{L/D}}{W}-\sin\gamma
\eq
that
restricts the admissible flight path angles in terms of the available
relative thrust.\RSfootnote{TWgamma01.m,singammatest.m}
\biglf
Once $V_{L/D}$ is defined, some texts, e.g.
\RScite{vinh:1980-1,hull:2007-1} introduce ``dimensionless speeds''
as ratios $V/V_{L/D}$, which are connected to our approach by
$$
\dfrac{V^2}{V^2_{L/D}}=\dfrac{R}{R_{L/D}}=R\dfrac{\sqrt{C_{D_0}}}{\sqrt{K}}
\cos \gamma=C_L\dfrac{\sqrt{C_{D_0}}}{\sqrt{K}},
$$
but the theory of quasi-steady flight
gets considerably simpler when using $R$.
\biglf
A {\em flight plan}
in the sense of this paper consists
of a function $h(x)$ of altitude $h$ over distance $x$. To turn it into a {\em
  trajectory}, an additional assignment
of speed or time along the flight plan is necessary.  
The {\em flight path angle} $\gamma$ is determined by 
$\tan \gamma(x)=\frac{dh(x)}{dx} =h'(x)$
independent of the speed assignment.
In view of our reduction of quasi-steady flight
to the $R$ variable, we shall consider 
assignments of $R$ instead of speed $V$
or time $t$ along the flight plan.
This way we split the calculation of optimal trajectories into two steps:
the determination of a speed or time assignment for each given flight plan, and
the variation of flight plans with given speed assignments.
\biglf
To deal with flight plans in terms of $x$, we go over to
the differential equation
\bql{eqWODE}
\begin{array}{rcl}
  \dfrac{dW}{dx}
  &=&\dfrac{\dot W}{\dot x}=\dfrac{-CT}{V\cos\gamma}\\
  &=&
-\sqrt{W}C\dfrac{T}{W}\dfrac{\sqrt{W}}{V\cos\gamma}\\
  &=&
-\sqrt{W}\dfrac{C}{\cos\gamma}\left(C_{D_0}R+
\dfrac{K\cos^2 \gamma}{R}+\sin \gamma\right)
\dfrac{\sqrt{\rho S}}{\sqrt{2R}}
\end{array} 
\eq
that, by substitution of $Z:=2\sqrt{W}$, turns into
a plain integration of the integrand
\bql{eqODE}
\dfrac{dZ}{dx}=-\dfrac{C}{\cos\gamma}\left(C_{D_0}R+
\dfrac{K\cos^2 \gamma}{R}+\sin \gamma\right)
\dfrac{\sqrt{\rho S}}{\sqrt{2R}}.
\eq
Solving this single ordinary differential equation
yields the weight along the flight plan,
and then speed and thrust follow via \eref{eqRconv} and \eref{eqTWfull}.
If needed, time $t(x)$ can be obtained by a parallel integration of
$1/V$ over $x$. 
We shall use \eref{eqODE} in various numerical examples,
once we have a strategy $R(h, \gamma)$ for choosing $R$.
A first case is \eref{eqRLD}, allowing to calculate
for any given flight plan $h(x)$ a speed assignment that realizes the
maximization of $L/D$ along the flight.
\biglf
But the right-hand side of \eref{eqODE}
also allows to calculate optimal flight plans for a given
strategy for $R(h, \gamma)$.
Indeed, if the right-hand side is written in terms
of $h(x)$ and $\gamma(x)=\arctan(h'(x))$, the minimization
of the integral leads to a variational problem that
has a second-order Euler-Lagrange equation whose solutions $h(x)$
minimize the integral, i.e. the overall fuel consumption.
We shall come back to this in Sections \RSref{SecVP} and \RSref{SecPS}.
%****************************************************************
\section{Range Maximization}\RSlabel{SecRMaxHF}
To maximize the range for a given amount of fuel or a prescibed weight loss,
one should introduce $W$ as the independent variable.
For a flight from position $x_0$ to $x_1$ with weight $W_0$ decreasing to $W_1$,
the distance covered is
$$
\int_{W_0}^{W_1}\dfrac{dx}{dW}dw.
$$
The integrand is
$$
%\begin{array}{rcl}
  \dfrac{dx}{dW}= \dfrac{\dot x}{\dot W}\\[0.3cm]
  =\dfrac{V\cos\gamma}{-CT}
%\end{array}
$$
and should be maximized. 
Before we do this optimization in general,
we consider a standard argument in the literature
\RScite{peckham:1974-1,vinh:1980-1,stengel:2004-1}
for the special case of {\em horizontal} flight.
There, 
$$
\begin{array}{rcl}
  \dfrac{dx}{dW}
  &=&
  %\dfrac{V}{-CD}\\
  %&=&-\dfrac{L}{D}\dfrac{V}{CW}\\
  %&=&-\dfrac{C_L}{C_D}
  %\dfrac{\sqrt{2W}}{\sqrt{C_L\rho S}CW}\\
  %&=&
  -\dfrac{\sqrt{C_L}}{C_D}
  \dfrac{\sqrt{2}}{C\sqrt{\rho SW}}
\end{array}
$$
leads to the conclusion that $\frac{\sqrt{C_L}}{C_D}$ is to be
maximized at each instant of an optimal horizontal quasi-steady flight.
Applying this to the drag polar \eref{eqdragpolar} yields
$$
{C_L}=\sqrt{\dfrac{C_{D_0}}{3K}},\; C_D=\frac{4}{3}C_{D_0}
$$
and a horizontal flight at a constant value of
\bql{eqR0}
R_0:=\sqrt{\dfrac{3K}{C_{D_0}}}=\sqrt{3}R_{L/D}
\eq
with a speed
\bql{eqV0}
V_0:=\sqrt[4]{3}\;V_{L/D}
\eq
that decreases with $\sqrt{W}$ like $V_{L/D}$. The same solution follows
when we express everything by $R$ via \eref{eqRconv} and minimize
the fuel consumption, i.e. the integrand in \eref{eqODE} with the major part
$$
%\begin{array}{rcl}
C_{D_0}\sqrt{R}+\dfrac{K}{R^{3/2}}
%\end{array} 
$$
over $R$.
This is a second strategy for determining $R$,
but it is restricted to horizontal flight, so
far.
\biglf
We now repeat this argument for general flight path angles,
and  use \eref{eqODE} to minimize
$$
%\begin{array}{rcl}
%  \dfrac{T}{V\cos\gamma}
%  &=&
%  \dfrac{
%  W\sqrt{\rho S}(C_{D_0}R+
%  \dfrac{K\cos^2 \gamma}{R}+\sin \gamma)}{\sqrt{2RW}\cos\gamma}\\
%  &=&\sqrt{W}
%  \dfrac{\sqrt{\rho S}}{\sqrt{2R}\cos\gamma}(C_{D_0}R+
%  \dfrac{K\cos^2 \gamma}{R}+\sin \gamma)\\
  %  &=&
  %\sqrt{W}
  %\dfrac{\sqrt{\rho S}}{\sqrt{2}\cos\gamma}(
  C_{D_0}\sqrt{R}+
  \dfrac{K\cos^2 \gamma}{R^{3/2}}+\dfrac{\sin \gamma}{\sqrt{R}}
  %)
%\end{array}
$$
over $R$
% that we already derived for \eref{eqODE}, 
with the solution\RSfootnote{Rgammasolve.mws}
\bql{eqRgamma}
R_\gamma:=\dfrac{1}{2C_{D_0}}
\left(\sin\gamma+\sqrt{\sin^2\gamma+12KC_{D_0}\cos^2\gamma}\right)
=R_0+\dfrac{\sin\gamma}{2C_{D_0}}+{\cal O}(\sin^2(\gamma)).
\eq
The other solution branch is always negative and unfeasible.
The solution could also be obtained in terms of $V$ or $T/W$,
but we can use our conversions \eref{eqRconv} and \eref{eqTWfull} to get
\bql{eqVgamma}
\begin{array}{rcl}
   \dfrac{T_\gamma}{W}
  &=&
  C_{D_0}R_\gamma+
  \dfrac{K\cos^2 \gamma}{R_\gamma}+\sin \gamma=:\tau(\gamma),\\[0.3cm]
  V_\gamma^2&=&\dfrac{2R_\gamma W}{\rho S}.
\end{array}
\eq
This gives an assignment of $R$ and speed $V$ for any given flight plan
by solving the ODE \eref{eqODE}. Like in \eref{eqRLD} and \eref{eqR0},
the resulting choice of $R$ depends only on $\gamma$ and the drag polar,
not on altitude and weight, which are built into $R$. 
\biglf
Before we vary these flight plans with $V_\gamma$ speed assignments
to get optimal trajectories, we add an illustration.
Figure \RSref{figTWRg}\RSfootnote{contours01.m}
shows the contours of the formula  \eref{eqTWfull} for $T/W$
plotted
in the $(\gamma,R)$ plane, for the values
$C_{D_0}=0.024$ and $K=0.073$ of the Standard Business Jet of
\RScite{hull:2007-1}.
The thick curve consists
of the points $(\gamma,R_\gamma)$ where $R$ is chosen optimally for given
$\gamma$. The leftmost
vertical line at $R_{L/D} =1.74$ hits all peaks of contour lines,
since $L/D$ maximization leads to the largest climb angle for a given
$T/W$, and all of these cases have the same $R$. The thick curve
meets this line at $T=0$, the engine-out situation,
where the optimal strategy is a glide at maximal $L/D$ ratio
with an angle of -4.78 degrees.
\biglf
The other vertical line
is at $R_0=\sqrt{3}R_{L/D}=3.02$ and hits the thick curve at $\gamma=0$, because
this is the well-known optimal choice of $R$ for horizontal flight.
The $T/W$ ratio for optimal horizontal flight is 0.0967, no matter
what the altitude, the weight, and the wing loading is. All of this is
coded into $R_0=3.02$, and the other 
variables can be read off \eref{eqRconv}.
\biglf
Each contour line resembles a special value of $T/W$ or a special
power setting chosen by the pilot. Then the points $(\gamma,R)$
on the contour describe the pilot's choice between climb angle and
speed (coded into $R$).
The maximal possible angle belongs to $R_{L/D}=1.74$, but this will
not be a good choice for range maximization. 
Of all the points on a given $T/W$ contour,
the intersection
of the contour with the thick $(\gamma,R_\gamma)$ curve
describes a special choice: at this $\gamma_{T/W}$,
the speed coded into $R_{\gamma_{T/W}}$
yields the optimum for range maximization. 
\begin{figure}[!tbp]
\begin{center}
\includegraphics[width=8cm,height=8cm]{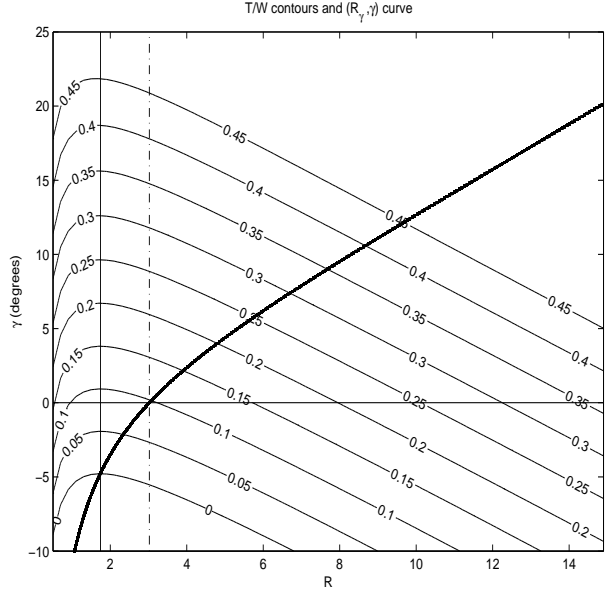}
\caption{Contours of $T/W$ over $(\gamma,R)$, with the
  optimality curve $(\gamma,R_\gamma)$
   \RSlabel{figTWRg}}
\end{center} 
\end{figure}
\biglf
The angle $\gamma_{T/W}$ can be calculated explicitly, because the
right-hand side of the equation
$$
\dfrac{T}{W}
  =
  C_{D_0}R_\gamma+
  \dfrac{K\cos^2 \gamma}{R_\gamma}+\sin \gamma
$$
  is the function $\tau(\gamma)$ from \eref{eqVgamma} 
  that can be inverted\RSfootnote{singammasolve.mws}
using MAPLE to yield
\bql{eqsingamW}
\sin(\gamma_{T/W})=\dfrac{2(T/W)
-\sqrt{(T/W)^2(1-12KC_{D_0})+64K^2C^2_{D_0}+16KC_{D_0}}}%
{2(1+4KC_{D_0})}.
\eq
If inserted into \eref{eqTWrestr}, the positive root
is infeasible. 
The above formula  can be applied to calculate a Continuous Descent
at nonzero idle thrust, or an optimal climb 
for a prescribed thrust policy as a function of altitude,
using the ODE \eref{eqODE} inserting $\gamma_{T/W}$ and $R_{\gamma_{T/W}}$.
We shall provide examples later. 
\biglf
If speed is prescribed, e.g. 250 knots indicated airspeed (KIAS)
below 10000 ft,
one has a prescribed $R$ and can use Figure \RSref{figTWRg} to read off
a $\gamma_R$ such that $R_{\gamma_R}=R$. The
only feasible solution\RSfootnote{Rminsolve.mws}  is
$$
6K\sin(\gamma_R)=R-\sqrt{R^2+36K^2-12KC_{D_0}R^2},
$$
and the thrust follows from  \eref{eqTWfull} again.
This will yield an optimal climb strategy under speed restriction,
solvable again via \eref{eqODE}. We shall come back to this in Section
\RSref{SecPS}.
%****************************************************************
  \section{Variational Problem}\RSlabel{SecVP}
  So far, we have determined the maximal-range instantaneous
  speed assignment for an arbitrary flight plan $h(x)$,
  given via $R_\gamma$ or $V_\gamma$ of \eref{eqRgamma} or \eref{eqVgamma}
for $\gamma=\arctan h'(x)$.
  If this speed does not violate
  restrictions, it is the best one for that flight plan. 
  But now we go a step further and vary the flight plans to find an optimal
  flight plan under all plans that allow the range-optimal instantaneous 
  speed assignment. 
\biglf
To this end, we insert $R_\gamma$ 
into the right-hand side of \eref{eqODE} 
to get a variational problem for the flight path $h(x)$. The
integrand for calculating $Z(x)=2\sqrt{W(x)}$ via
\bql{eqvarprob}
Z(x_1)-Z(x_0)=2\sqrt{W(x_1)}-2\sqrt{W(x_0)}=\int_{x_0}^{x_1}\dfrac{dZ}{dx}dx
\eq
is
\bql{eqODE2}
\dfrac{dZ}{dx}=-S\dfrac{C(h)\sqrt{\rho(h)}}{\sqrt{2}}
\dfrac{C_{D_0}R_\gamma+
  \dfrac{K\cos^2 \gamma}{R_\gamma}+\sin \gamma}{\sqrt{R_\gamma}\cos\gamma}
=:F(h)G(h'),
\eq
and the variational problem consists of finding $h(x)$ such that
the integral in \eref{eqvarprob} is minimized.
The integrand is a product of a function $F$ of $h$ and
a function $G$ of $h'$ via
$\gamma=\arctan(h')$. For such a variational problem, 
the Euler-Lagrange equation is 
\bql{eqEulerLag}
h''=\dfrac{F'(h)}{F(h)}\left( \dfrac{G(h')}{G''(h')}-
\dfrac{G'(h')}{G''(h')}h'\right)
\eq
by standard arguments of the Calculus of Variations,
and we need the corresponding complicated derivatives of $F$ and $G$.  
\biglf
The function $G$ is dependent only on the drag polar, not on
propulsion, and derivatives wrt. $h'$ can be
generated by symbolic computation, e.g. using MAPLE.
The function $F(h)$ is $C(h)\sqrt{\rho(h)}$ up to constants and depends
on propulsion
only via the altitude-dependency of the specific fuel consumption $C(h)$.
In simple models, e.g. \RScite{hull:2007-1} for turbofans and turbojets,
$C(h)$ is an exponential function of $h$, as well as the air density $\rho(h)$.
Then symbolic computation will work as well for the $h$-dependent part.
Using the code generation feature of MAPLE, one gets ready-to-use
expressions in MATLAB for solving the second-order ODE \eref{eqEulerLag}
for optimal flight plans $h(x)$,
without any detour via Optimal Control. 
\biglf
A closer inspection of the Euler-Lagrange equation 
for the variational problem shows that $F'(h)/F(h)$ is a constant   
if $C(h)$ and $\rho(h)$ have an exponential law, and then the
right-hand-side of the Euler-Lagrange equation \eref{eqEulerLag} 
is a pure equation in $h'$.
Since the equation is also autonomous, i.e. independent of $x$,
the solutions $h(x)$ in the $(x,h)$
plane can be shifted right-left and up-down.
\begin{figure}[!tbp]
\begin{center}
\includegraphics[width=6.0cm,height=6.0cm]{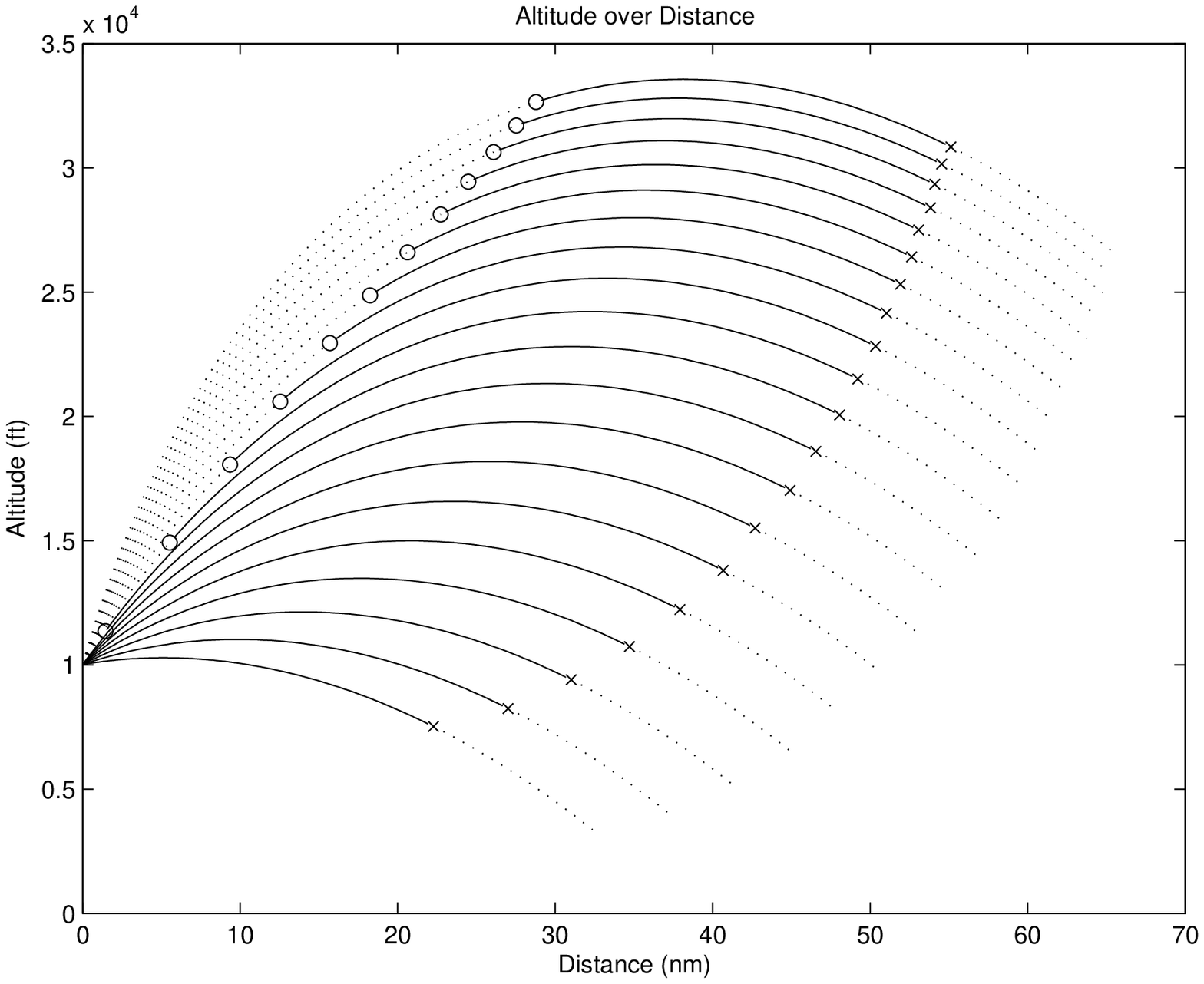}
\includegraphics[width=6.0cm,height=6.0cm]{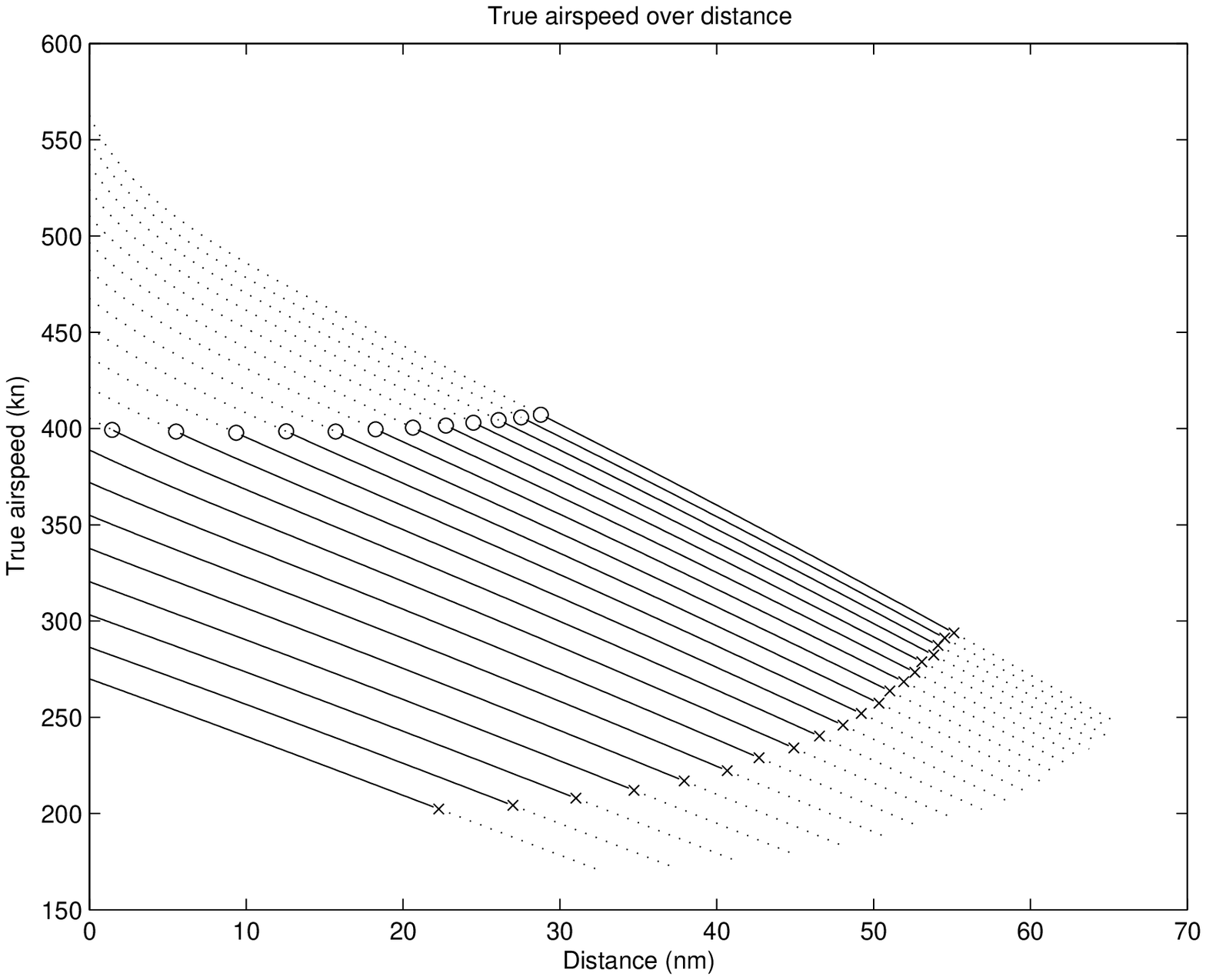}\\
\includegraphics[width=6.0cm,height=6.0cm]{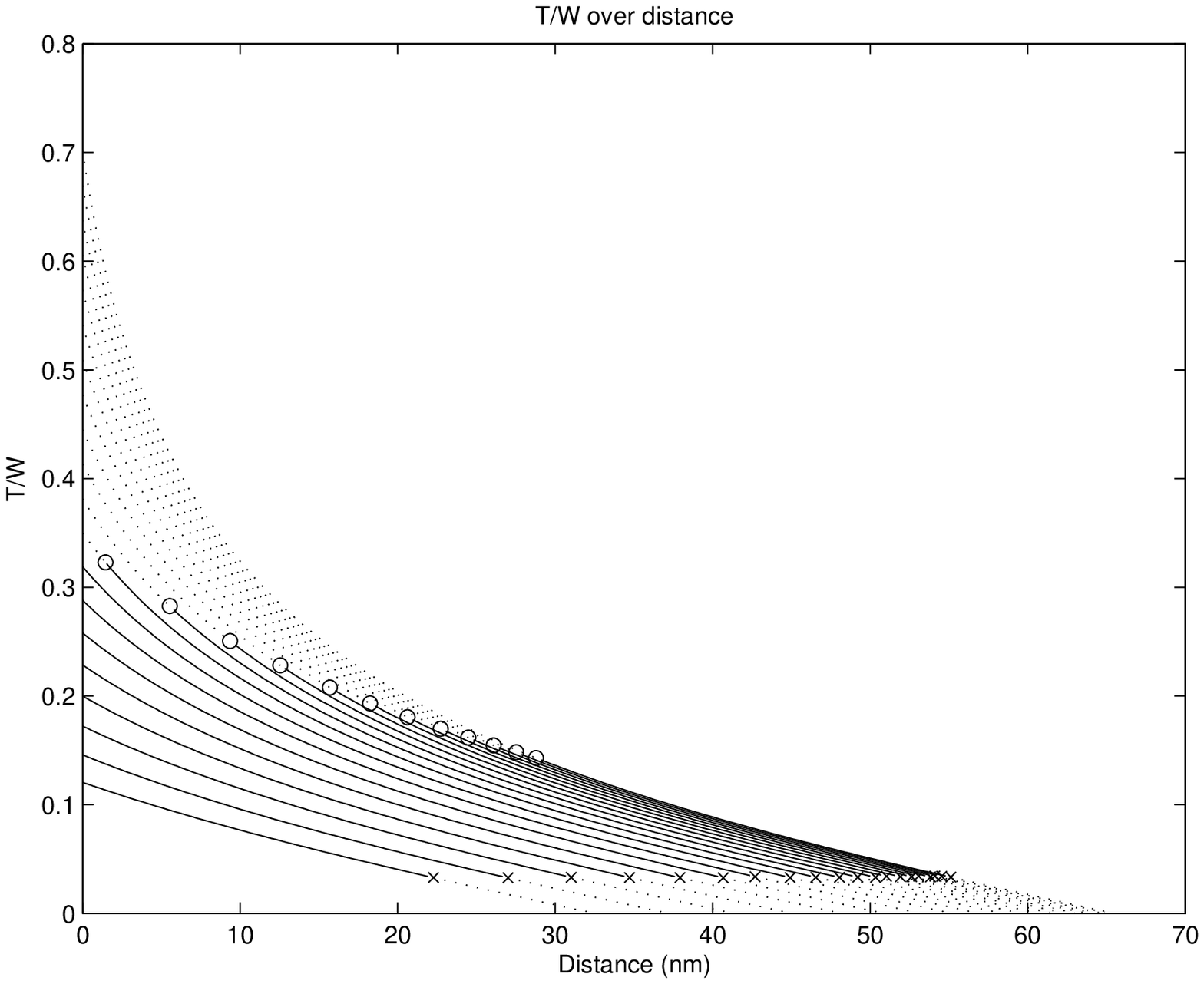}
\includegraphics[width=6.0cm,height=6.0cm]{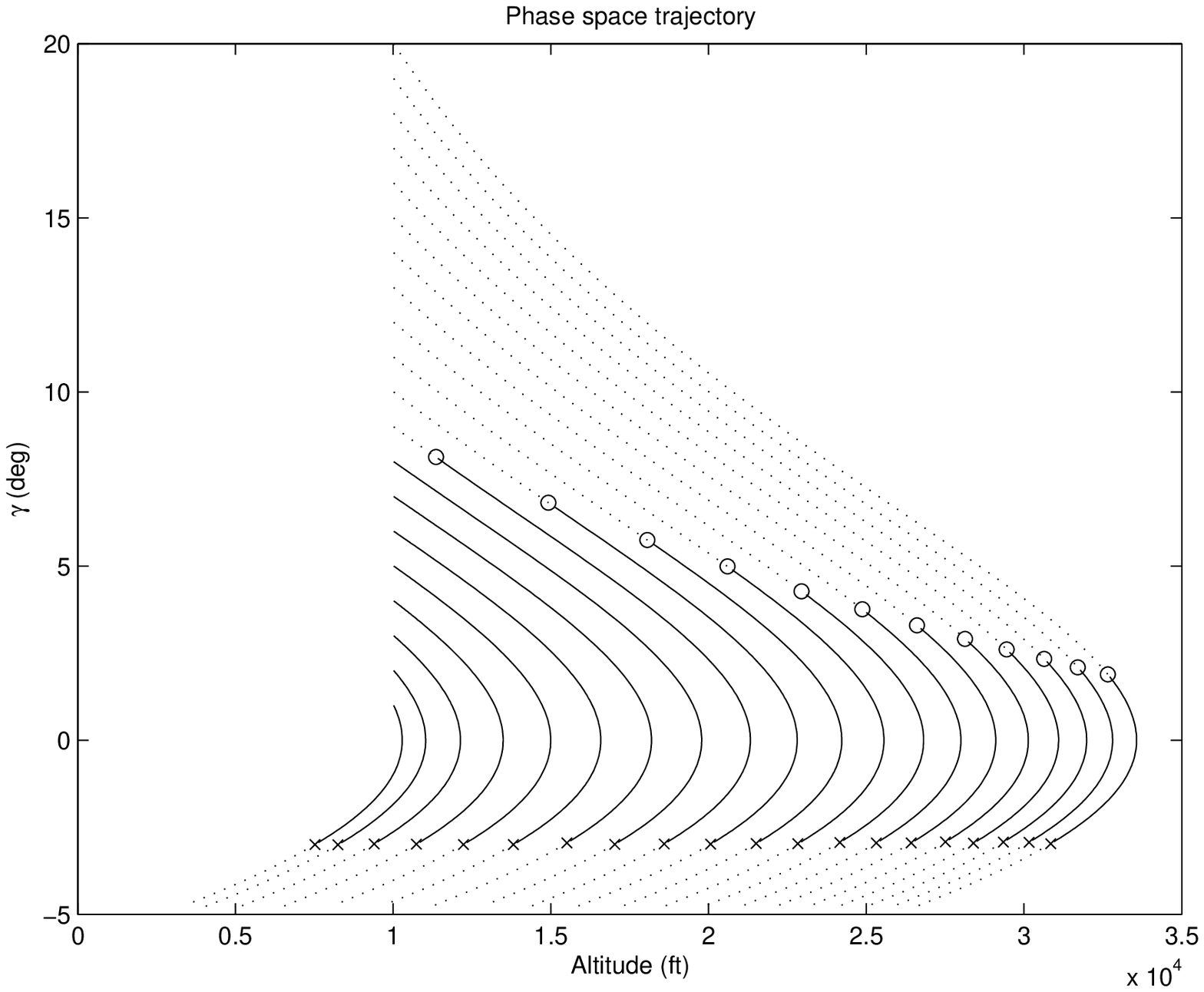}
\caption{Some unconstrained flight paths
  for the SBJ model in \RScite{hull:2007-1}, plots of 
  $(h,x),\;(V,x),\;(T/W,x),\;(\gamma,h)$
   \RSlabel{FigUnconFP}}
\end{center} 
\end{figure}
\biglf
Figure\RSref{FigUnconFP}\RSfootnote{testUCFP.m}
shows typical solutions of the
Euler-Lagrange equation for the Standard Business Jet (SBJ)
model from
\RScite{hull:2007-1}, starting at 10000 ft and ending at 3000 ft.
A closer inspection of the
differential equation reveals that the solutions are always concave
in the $(x,h)$ plane, and the
speed is always decreasing, see the two upper plots.
The lower left plot shows the $T_\gamma/W$ values
of \eref{eqVgamma}, and these may be
too large or too small to be admissible. Therefore all curves are 
dotted where the thrust restrictions are violated.
The lower right plot of Figure \RSref{FigUnconFP} 
visualizes this in phase space, where we
replaced $h'$ by $\gamma=\arctan h'$
for convenience.  The trajectories there
are traversed downwards, with decreasing $\gamma$,
and the extremum of $h$
to the right.
\biglf
This looks disappointing at first sight,
but we have to take the thrust
limits into account and view the variational problem as a
{\em constrained} one. Such problems have
the well-known property that solutions either
follow
the Euler-Lagrange equation or a
boundary defined by the restrictions.
In our case, only the solid curves
between the circles and the crosses
are 
solutions of the Euler-Lagrange equations that solve the
unconstrained variational problem. When a solution of
the variational problem hits a constraint, the Euler-Lagrange
ODE is not valid anymore, but one can use the constraint
to determine the solution. We shall do that in what follows,
and point out that optimal full flight plans will
follow the circles first, then
depart from the circles to a solid line, and
depart form the line at a cross
to follow the crosses from that point on. This
argument is qualitatively true, but needs a minor
modification due to the fact that the true
weight behaves slightly differently
when we consider a single trajectory, while 
Figure \RSref{FigUnconFP} shows multiple trajectories. 
%****************************************************************
\section{Constrained Range-Optimal Trajectories}\RSlabel{SecCROT}
We now check the solutions of the variational problem
when thrust restrictions are active.
These are partial
flight plans with range-optimal speed assignments as well
as the partial flight plans that do not violate restrictions,
being solutions of the Euler-Lagrange equation
\eref{eqEulerLag}. To calculate the thrust-restricted parts,
we assume thrust being given as a function of
altitude, either as maximal admissible
continuous thrust or as
idle thrust. Inserting the current weight $W$,
we use \eref{eqsingamW} to calculate the flight path angle
$\gamma_{T/W}$ that yields the range-optimal
$R_{\gamma_{T/W}}$ assignment via  \eref{eqRgamma}.
Then an ODE system for $h$ and $W$ is set up using
\eref{eqWODE} and $h'(x)=\tan \gamma_{T/W}$.
\biglf
Doing this for maximal admissible continuous thrust
yields range-optimal climb/cruise trajectories, while
inserting idle thrust yields range-optimal Continuous
Descent trajectories.  Between these two parts of a
range-optimal flight, there must be a {\em transition}
from maximal admissible continuous thrust to idle thrust,
and this transition must follow a
solution of the Euler-Lagrange equation. In
terms of Figure \RSref{FigUnconFP}, the climb/cruise path
reaches a circle, then follows one of the curves up to the cross
marking idle thrust,
and then a Continuous Descent trajectory follows. The
{\em Top of Descent} point is reached in the transition part.  
\biglf
Starting at a given altitude and weight, the speed and the initial
flight path angle are determined.\RSfootnote{Vgamma.m}
Because the range-optimal speed
$V_\gamma$ usually comes out to be well above 250 KIAS
at low altitudes, we start our range-optimal trajectories at 10000
ft, and for the following plots we used a fixed starting weight
at that altitude.
\begin{figure}[!tbp]
\begin{center}
\includegraphics[width=6.0cm,height=6.0cm]{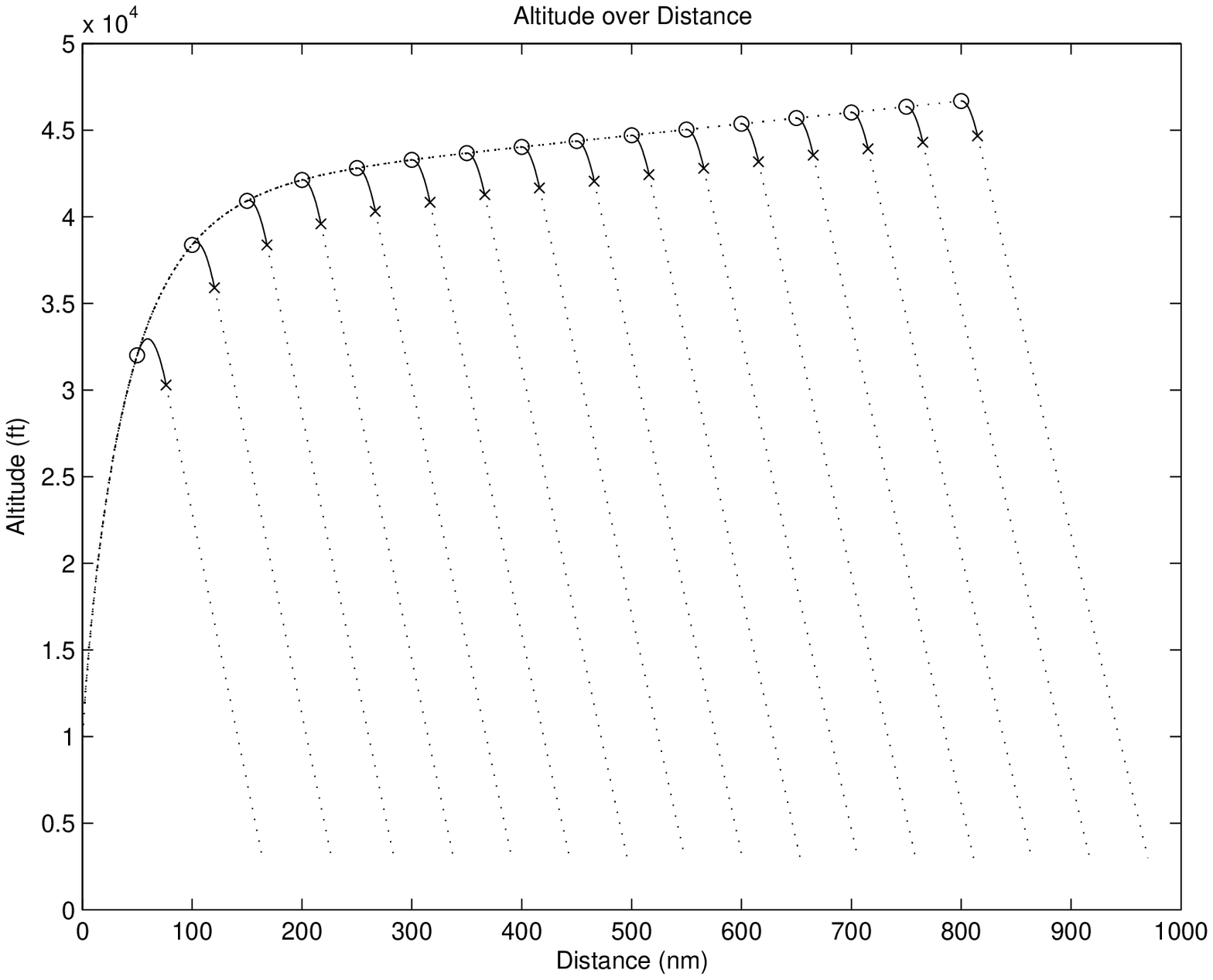}
\includegraphics[width=6.0cm,height=6.0cm]{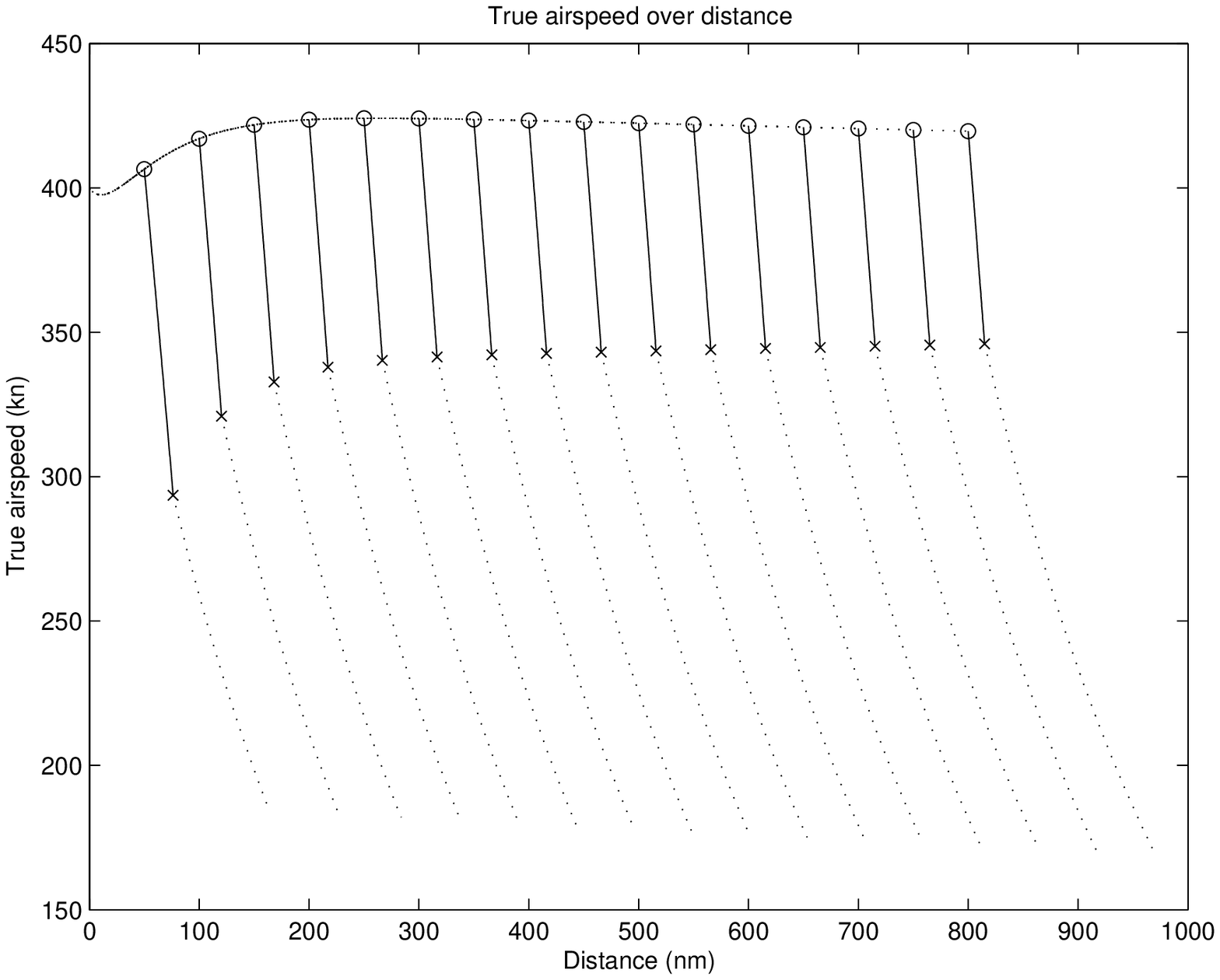}\\
\includegraphics[width=6.0cm,height=6.0cm]{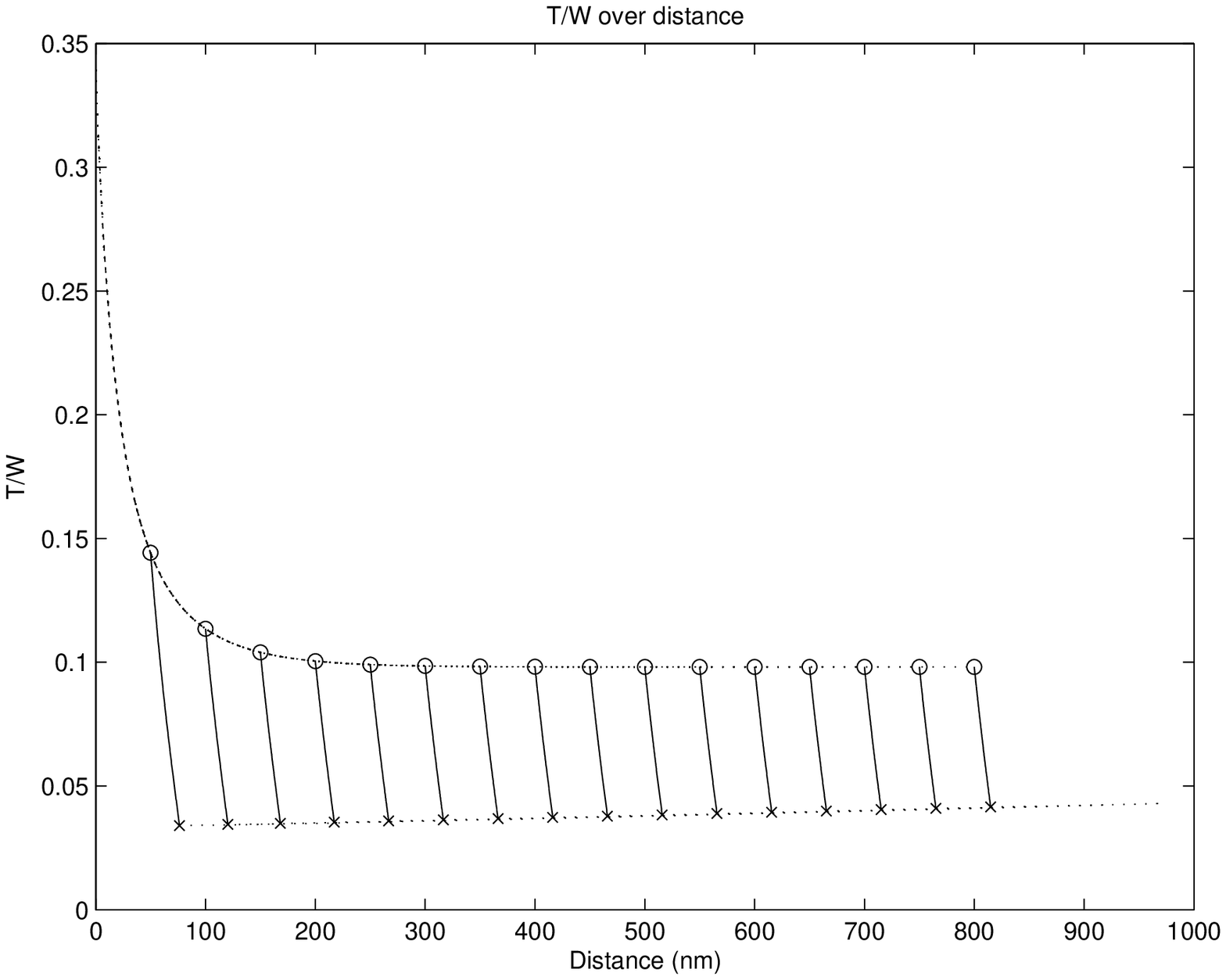}
\includegraphics[width=6.0cm,height=6.0cm]{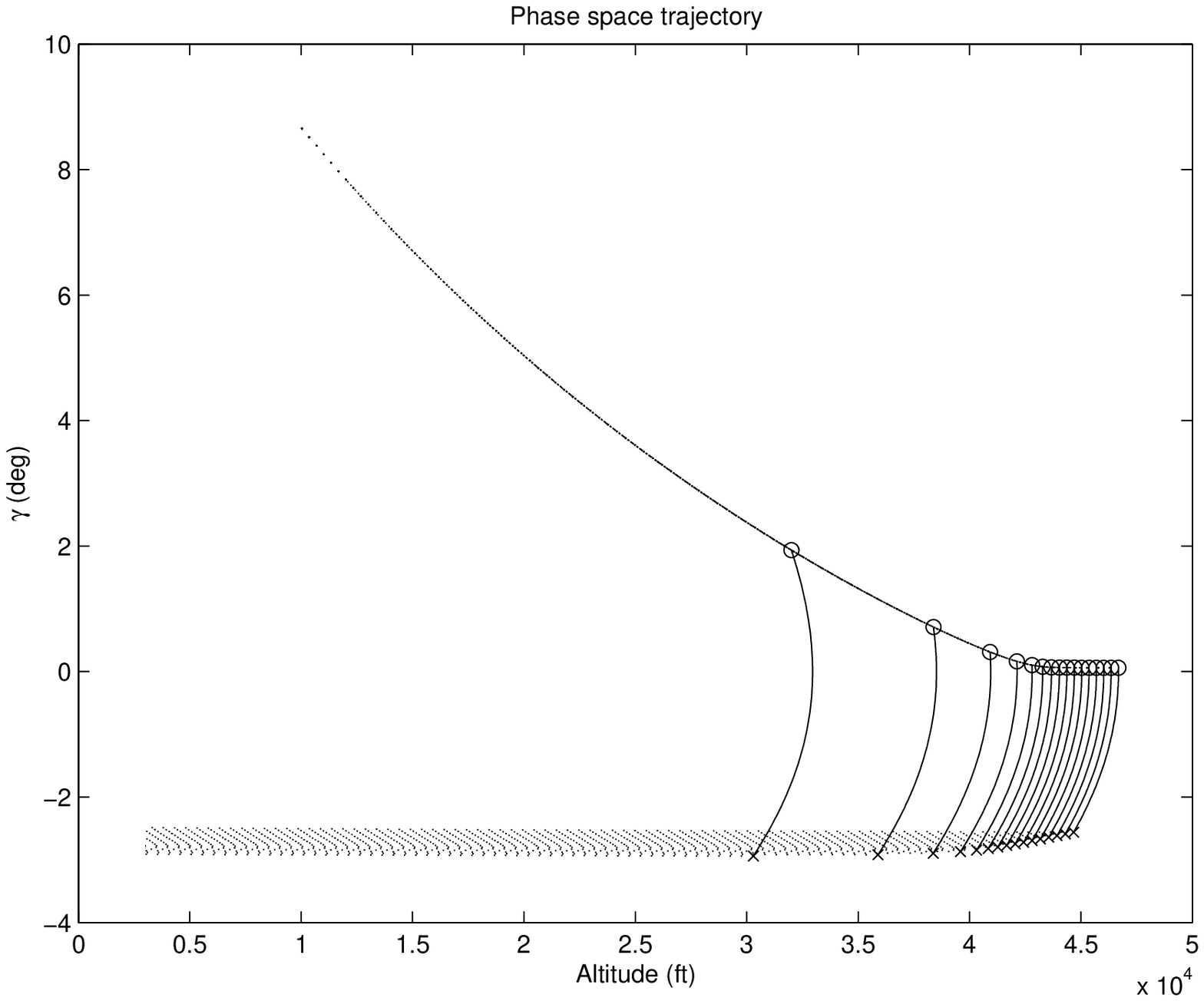}
\caption{Some optimal three-piece flight paths
  for the SBJ model in \RScite{hull:2007-1}, plots of 
  $(h,x),\;(V,x),\;(T/W,x),\;(\gamma,h)$
   \RSlabel{Fig3pFP}}
\end{center} 
\end{figure}
\biglf
Figure \RSref{Fig3pFP}\RSfootnote{threepiece.m}
shows range-optimal trajectories with the three
parts described above, for the Small Business Jet of
\RScite{hull:2007-1}. The climb/cruise part, using a maximal continuous thrust
power setting of 0.98,  is stopped at distances
from 50 to 800 nm in steps of 50 nm to produce the different
trajectories. When a transition is started,
the final $\gamma$ of the
climb is used to calculate an unconstrained solution of the
Euler-Lagrange equation that performs a smooth transition to
the Continuous Descent part at idle thrust. Along the Euler-Lagrange transition,
the decreasing $T/W$ values are monitored,
and the Continuous Descent is started
when $T_{idle}/W$ is reached. The full range-optimal flight paths are in the
top left plot, while the top right shows the true airspeed
and
the bottom left shows the $T/W$ values along the flight paths.  
The final plot is in phase space. One can compare with
Figure \RSref{FigUnconFP}, but there the
total flight distances are much smaller.
To arrive at a certain destination distance and altitude, the starting point
of the transition has to be adjusted. 
\biglf
A close-up of one of the transitions is in
Figure\RSfootnote{threepiecesingle.m} \RSref{Fig3pFPsingle},
namely the one where the transition is started at 400 nm. The transition
takes about 15 nm, and the right-hand plot shows what the pilot
should do for a range-optimal flight: decrease thrust from
maximal continuous thrust to idle thrust slowly and roughly linearly,
using about 15 nm.
At high altitudes, the top-of-descent point is reached very
shortly after the transition is started, see the phase space plot
in Figure \RSref{Fig3pFP}. 
\begin{figure}[!tbp]
\begin{center}
\includegraphics[width=6.0cm,height=6.0cm]{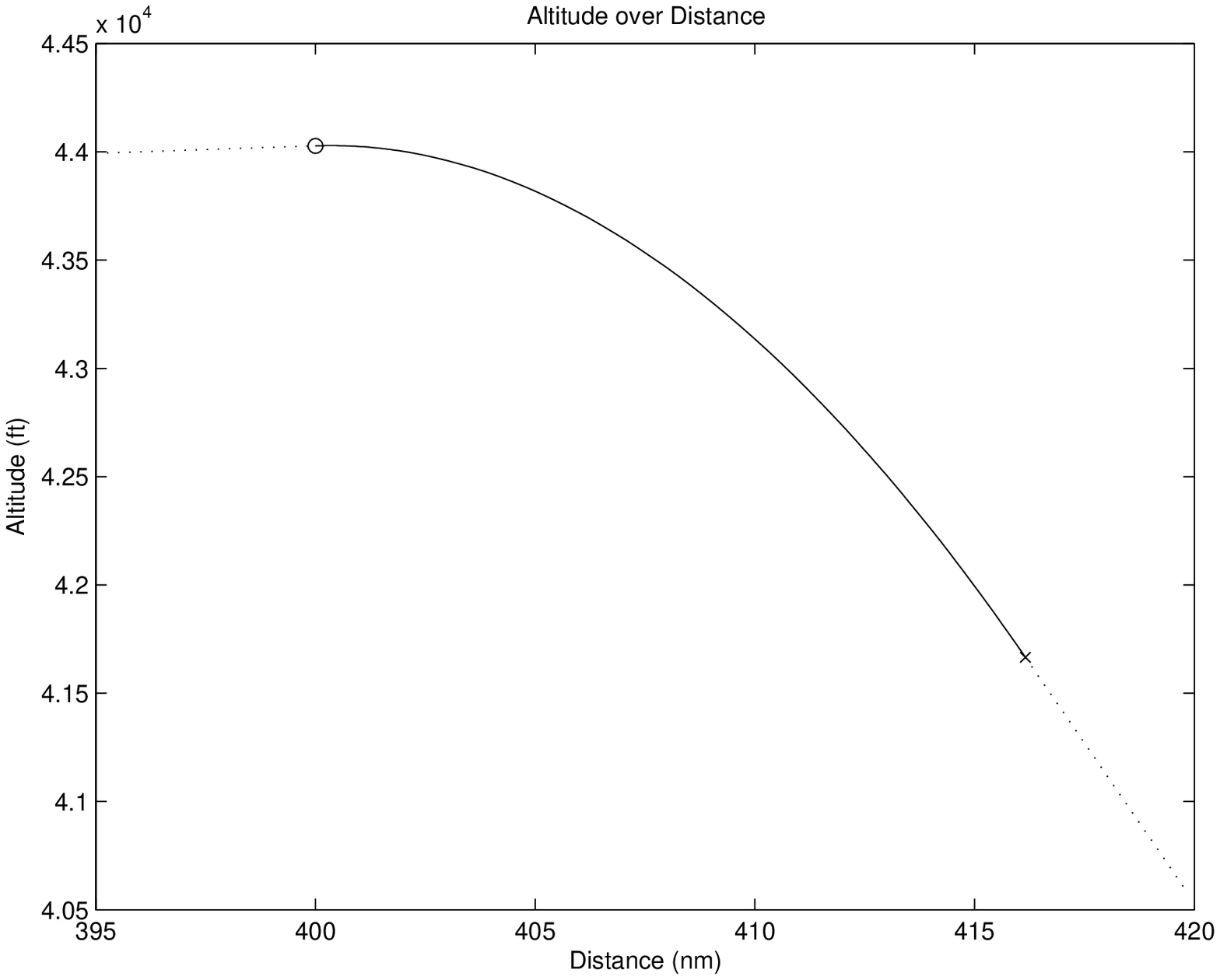}
\includegraphics[width=6.0cm,height=6.0cm]{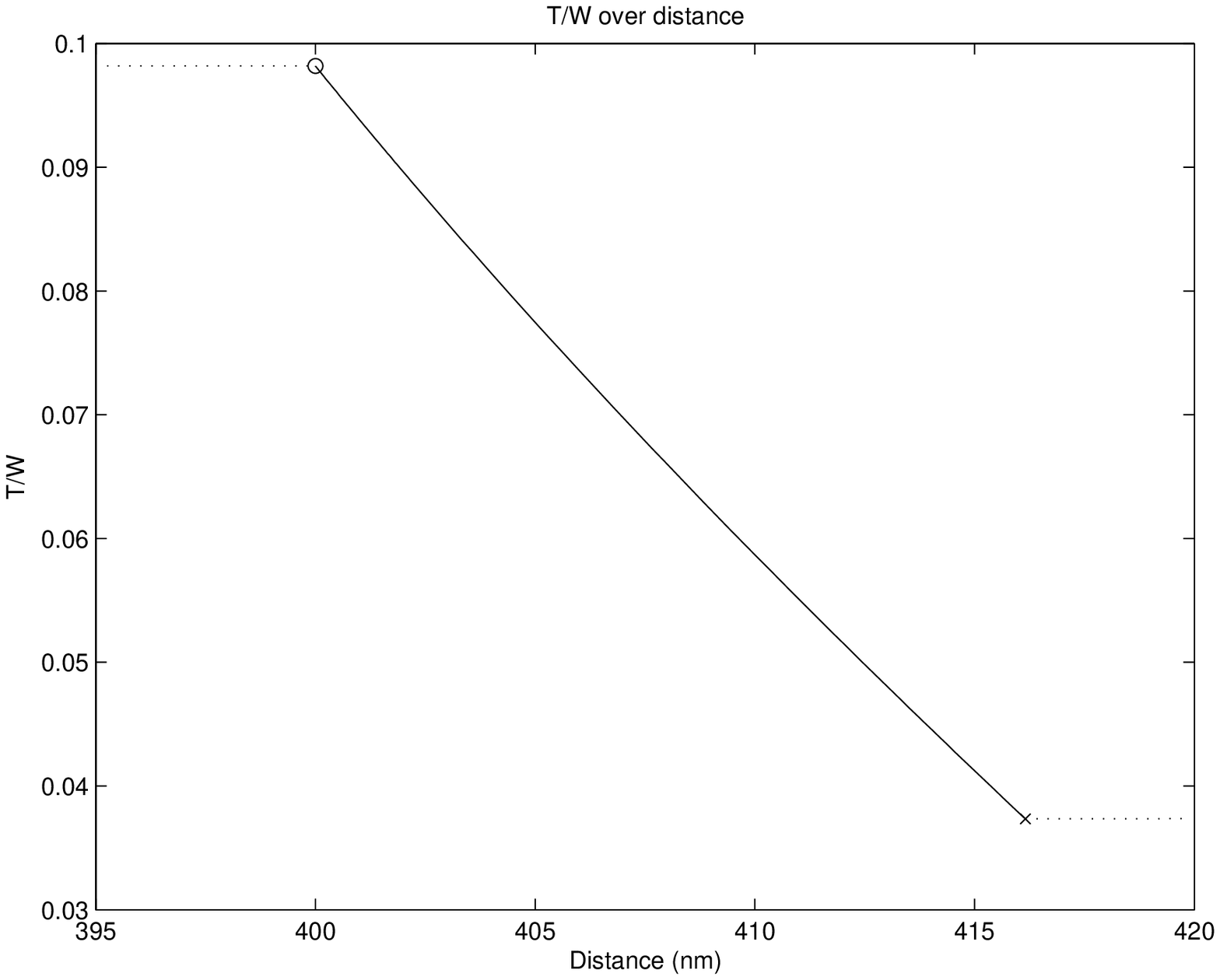}
\caption{Transition part of an optimal three-piece flight path
  for the SBJ model in \RScite{hull:2007-1}, plots of 
  $(h,x),(T/W,x)$
   \RSlabel{Fig3pFPsingle}}
\end{center} 
\end{figure}
\biglf
This means that range-optimal long-distance flights above 10000 ft
have necessarily three sections:
\begin{enumerate}
\item a climb/cruise at maximal continuous admissible thrust,
%  proceeding in phase space from left to right along the upper red/blue
%  boundary, 
\item a transition following a solution of the Euler-Lagrange equation,
%  downward along a blue curve in phase space,
\item and a continuous descent at idle thrust.
%  , along the lower red/blue
%    boundary, going from right to left. 
\end{enumerate}
\biglf
Figure \RSref{Fig3pFP} shows that for high altitude the flight path angle
tends to be constant. To analyze this effect, we go over to a
single differential equation for $\gamma$ that has a stationary solution.
% \section{Direct Try}\RSlabel{SecDT}
We insert the prescribed thrust $T(h)$ 
into \eref{eqVgamma} to get\RSfootnote{taudef.mws}
$$
\dfrac{T(h)}{W}=\tau(\gamma)=2
\dfrac{s^2+s\sqrt{s^2+12KC_{D_0}c^2}+4KC_{D_0}c^2}{s+\sqrt{s^2+12KC_{D_0}c^2}}
$$
with $\tau$ from \eref{eqVgamma} being
the inverse function of \eref{eqsingamW} in terms of
$s=\sin\gamma$ and $c=\cos\gamma$. 
The idea now is to get rid of $W$ by taking $h$-derivatives and 
$$
\begin{array}{rcl}
  T'(h)&=& \dfrac{d}{dh}(W\tau(\gamma))
%  &=&
%  \dfrac{dW}{dh}\tau(\gamma) +W\tau'(\gamma)\dfrac{d\gamma}{dh}\\[0.3cm]
  =
  \dfrac{dW}{dh}\tau(\gamma) +\dfrac{T(h)}{\tau(\gamma)}
  \tau'(\gamma)\dfrac{d\gamma}{dh},\\[0.3cm]
%\end{array}
%$$
%$$
%\begin{array}{rcl}
  V^2
  &=&
  \dfrac{2R(\gamma)W}{\rho(h)S}
  =
  \dfrac{2R(\gamma)T(h)}{\tau(\gamma)\rho(h)S}\\[0.3cm]
  \dfrac{dW}{dh}
  &=& \dfrac{\dot W}{\dot h}
  = \dfrac{-C(h)T(h)}{V\sin\gamma}
  =
  -\dfrac{C(h)\sqrt{T(h)\tau(\gamma)\rho(h)S}}{\sin\gamma\sqrt{2R(\gamma)}}
\end{array}
$$
to arrive at the single differential equation
$$
%\begin{array}{rcl}
  \dfrac{T'(h)}{T(h)}
%  &=&  -\dfrac{C(h)\sqrt{\tau(\gamma)\rho(h)S}}{\sin\gamma\sqrt{2T(h)R(\gamma)}}\tau(\gamma)
%  +\dfrac{\tau'(\gamma)}{\tau(\gamma)}\dfrac{d\gamma}{dh}\\[0.3cm]
%  &=&  -\dfrac{C(h)\sqrt{\rho(h)}\tau^{3/2}(\gamma)\sqrt{S}}%
%  {\sqrt{T(h)}\sin\gamma \sqrt{R(\gamma)}\sqrt{2}}
%  +\dfrac{\tau'(\gamma)}{\tau(\gamma)}\dfrac{d\gamma}{dh}%\\[0.3cm]
  %
  =  -\dfrac{C(h)\sqrt{\rho(h)}}{\sqrt{T(h)}}
  \dfrac{\tau^{3/2}(\gamma)}{\sin\gamma \sqrt{R(\gamma)}}
  \dfrac{\sqrt{S}}{\sqrt{2}}%
  +\dfrac{\tau'(\gamma)}{\tau(\gamma)}\dfrac{d\gamma}{dh}\\[0.3cm]
%\end{array}
  $$
  governing range-optimal climb/cruise at prescribed thrust. 
  If a constant $\gamma$ would solve this ODE,
  the equation
$$
  \dfrac{-T'(h)}{\sqrt{T(h)}C(h)\sqrt{\rho(h)}}
  =\dfrac{\tau^{3/2}(\gamma)}{\sin\gamma \sqrt{R(\gamma)}}
  \dfrac{\sqrt{S}}{\sqrt{2}}
  $$
  must hold over a certain range of $h$.
  But for the turbojet/turbofan propulsion models
  of \RScite{hull:2007-1}
  and the exponential air density model \eref{eqairdens},
  all parts of the left-hand side are certain powers
  of $\rho(h)$ that finally cancel out, letting the 
  left-hand side be a constant that only depends on the power setting.
  The right-hand side has a singularity for $\gamma=0$, and 
  there always
  is a small fixed positive angle $\gamma$ solving the above equation.
  The ODE solution tends to this 
  for increasing $h$, explaining the constant 
  final climb angle in Figure \RSref{Fig3pFP}
  for prescribed thrust. The left-hand side seems to be a
  crucial parameter for
  propulsion design, relating consumption to thrust and altitude
  for both turbojets and turbofans.%
  \RSfootnote{test01.m, testhpart.m, Direkt02.mws, Direkt04.mws}
%****************************************************************
\section{Prescribed Speed}\RSlabel{SecPS}
To deal with the usual speed restriction below 10000 ft,
we have to abandon the above
scenario, because we cannot minimize
fuel consumption with respect to speed
anymore. If the speed is given (in terms of $R$), equation
\eref{eqTWfull} still has one degree of freedom,
connecting $T/W$ to the flight path angle $\gamma$,
and we have to solve for a
range-optimal climb strategy in a different way now, going directly to a
variational problem.
\biglf
Implementing a 250 KIAS 
restriction
including conversion to true airspeed,
we have an altitude-dependent prescribed speed $V_F(h)$. 
Since air density $\rho$
is also $h$-dependent, so is the dynamic pressure
$\bar q(h)=\frac{1}{2}\rho(h) V_F(h)^2$ and
the variable $U(h)=\bar q(h)S$ connecting the
pressure ratio $R$ to the weight
$W$ via 
$$
R(h,W)=\dfrac{\bar q(h)}{\frac{W}{S}}
=\dfrac{\rho(h)V_F^2(h)S}{2W}=\dfrac{U(h)}{W}.
$$
Then \eref{eqTWfull} yields
$$
\begin{array}{rcl}
%  \dfrac{T}{W}
%  C_{D_0}R+\dfrac{K\cos^2\gamma}{R} +\sin\gamma\\[0.4cm]
  T=
%  &=&
%  C_{D_0}RW+\dfrac{KW^2\cos^2\gamma}{RW} +W\sin\gamma\\[0.4cm]
%  &=&
  C_{D_0}U(h)+\dfrac{KW^2\cos^2\gamma}{U(h)} +W\sin\gamma.
\end{array} 
%T=WC_{D_0}\bar q(h)S+\dfrac{W^2K\cos^2\gamma}{\bar q(h)S}+W\sin\gamma.
$$
Inserting into the fuel consumption integrand, we get
$$
\begin{array}{rcl}
  -{W'(x)}
  &=&
%  \dfrac{C(h)T}{V_F(h)\cos\gamma}\\[0.4cm]
%  &=&
%  \dfrac{C(h)}{V_F(h)\cos\gamma}\left(
% C_{D_0}U(h)+\dfrac{KW^2\cos^2\gamma}{U(h)} +W\sin\gamma \right)\\[0.4cm] 
%  &=&
  \dfrac{C(h)}{V_F(h)}\left(
 \dfrac{C_{D_0}U(h)}{\cos\gamma}+\dfrac{KW^2\cos\gamma}{U(h)} +W\tan\gamma \right)\\[0.4cm] 
 \end{array}
$$
which is  a Lagrangian $L(W,h,h')=L(W,u,v)$ and leads to a variational problem
with an Euler-Lagrange equation. In contrast to Section
\RSref{SecVP}, the weight is not eliminated, but we simply keep
it in the Lagrangian.
We only need the Lagrangian for $h$ up to 10000 ft,
and then we can fit each $h$-dependent part with good accuracy by a
low-degree polynomial in $h$.The $\gamma$- or $h'=\tan\gamma$- dependent
parts can be differentiated symbolically, as well as the
polynomial approximations to the $h$-dependent parts. We get the ODE system
$$
\begin{array}{rcl}
  h'&=& v\\
  v'&=&\dfrac{1}{L_{vv}(W,u,v)}(L_{u}(W,u,v)-L_{vu}(W,u,v)v)\\
  W'&=&-L(W,u,v)
\end{array}
$$
for the Euler-Lagrange
flight paths, under suitable initial or boundary value
conditions, and we can roughly repeat Section
\RSref{SecVP} for the new variational problem.
Above, the subscripts denote the partial derivatives. 
\begin{figure}[!]
\begin{center}
\includegraphics[width=4cm,height=4cm]{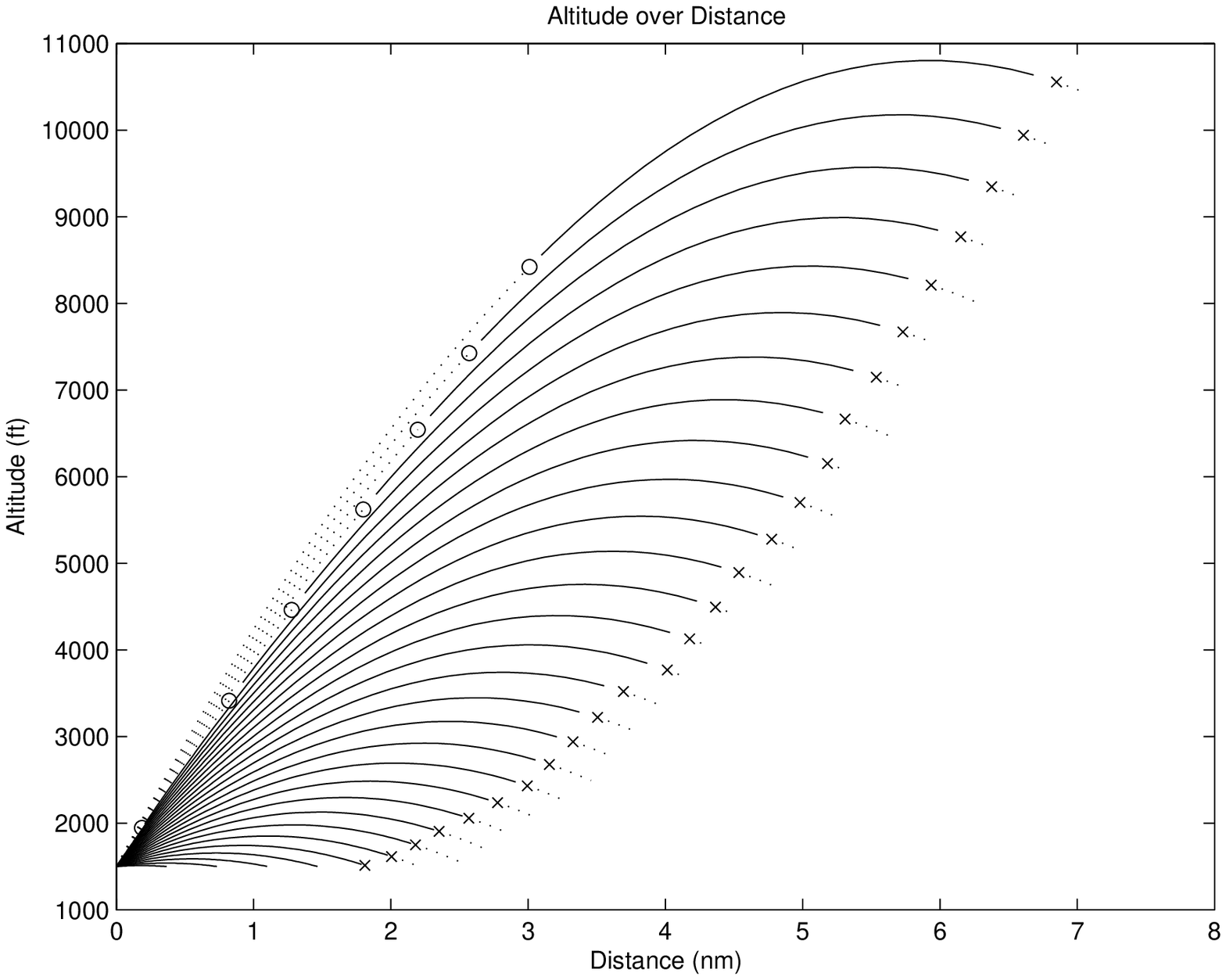}
\includegraphics[width=4cm,height=4cm]{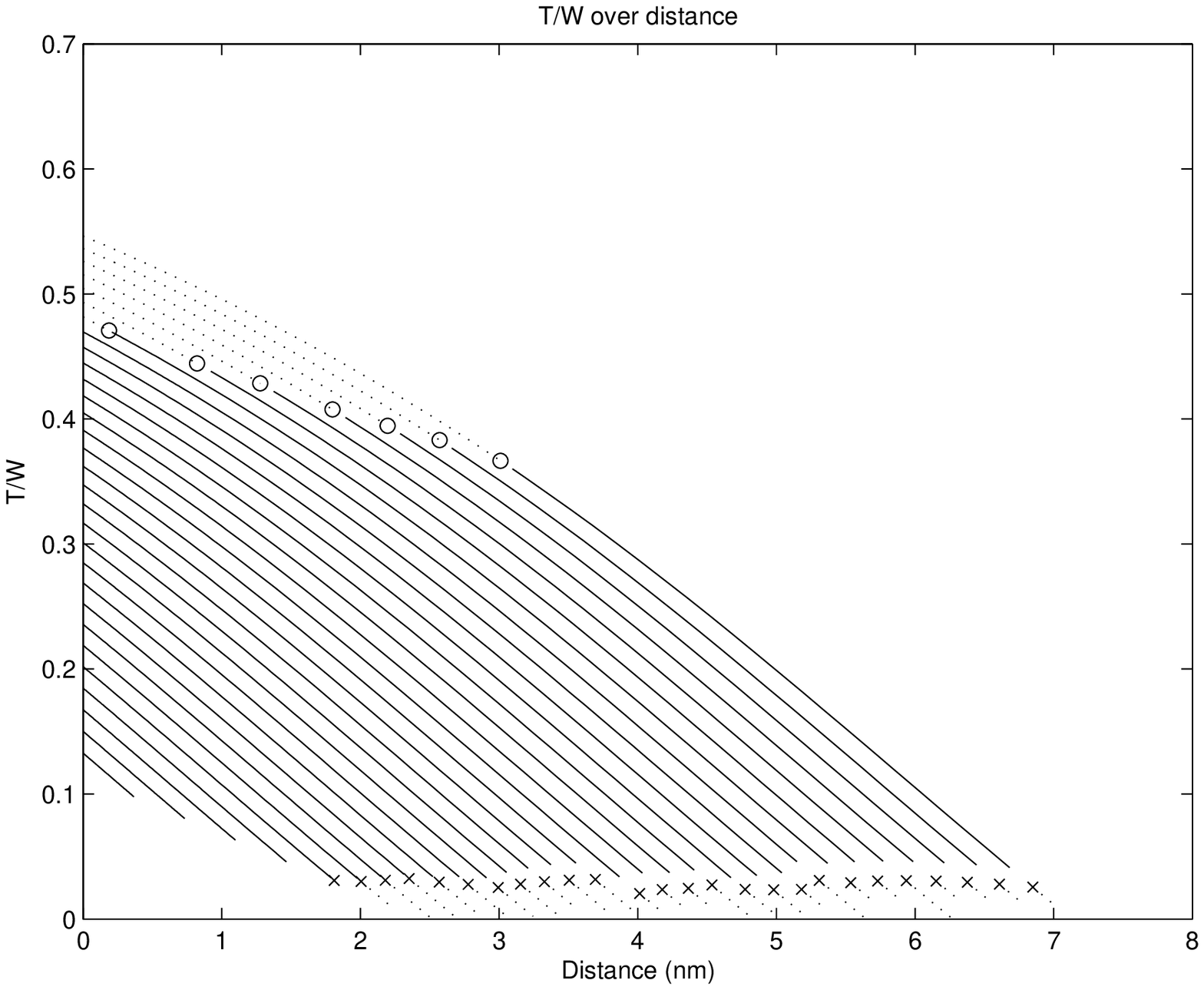}
\includegraphics[width=4cm,height=4cm]{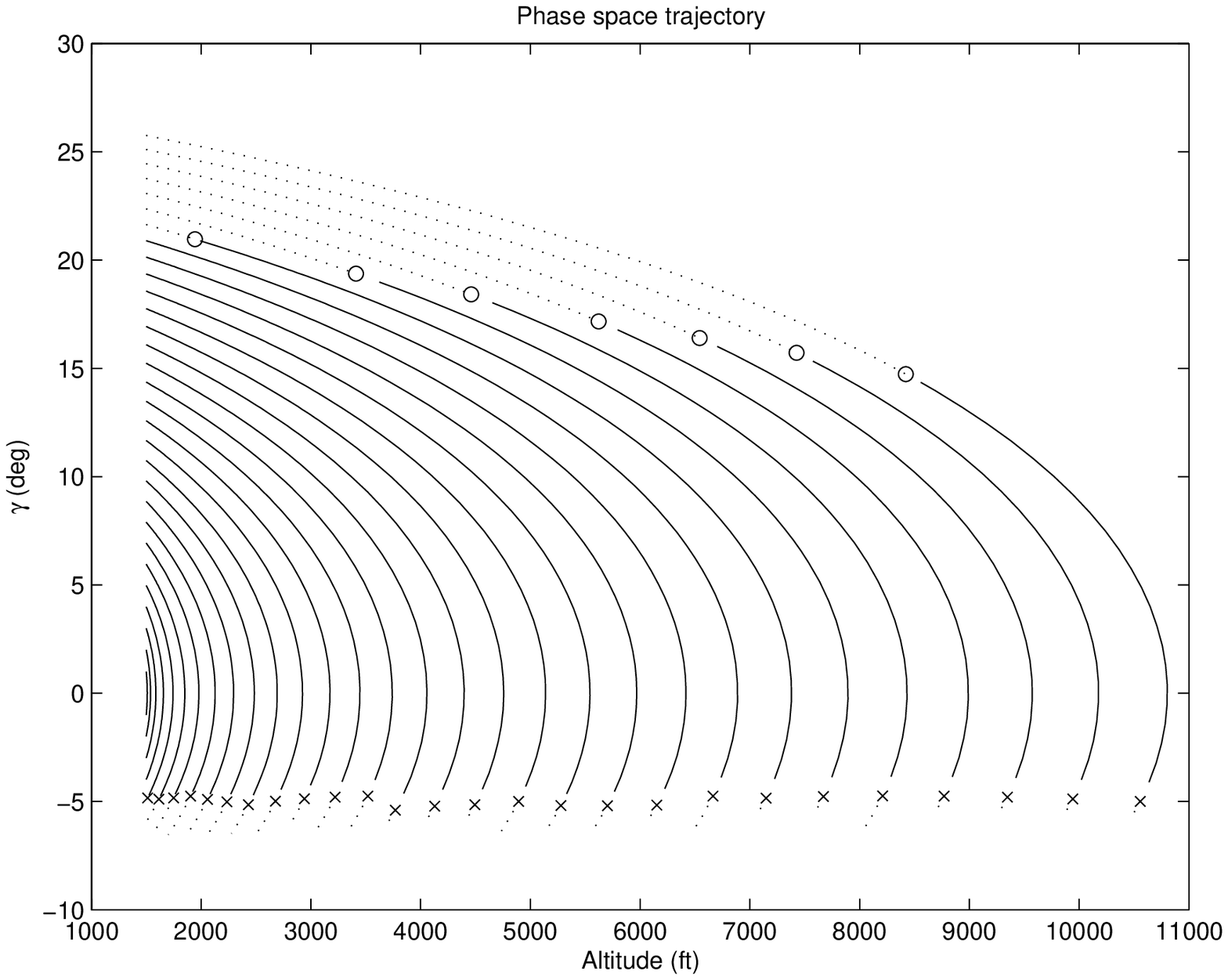}
\caption{Unconstrained range-optimal solutions
at 250 KIAS,  flight path, $T/W$ ratio, and phase space 
   \RSlabel{FigFVUCFP}}
\end{center} 
\end{figure}
\biglf
For the aircraft model in \RScite{hull:2007-1}
started at 1500 ft with maximal weight, we get Figure
\RSref{FigFVUCFP}\RSfootnote{FVUCFP.m}
showing range-optimal
unconstrained trajectories for flight at 250 KIAS at low altitudes.
Like in the previous figures, the dotted parts violate thrust restrictions.
For a long-range flight, the trajectory
reaching $\gamma=0$ exactly at 10000 ft should be selected,
and but it needs excessive thrust at the beginning.
\biglf
Therefore the upper
thrust limit for the variational problem has to be accounted for,
and range-optimal trajectories for a 250 KIAS climb will
consist of two pieces: the first with maximal admissible thrust,
and the second as a transition satisfying the Euler-Lagrange equation
for the optimal speed-restricted case.
Because the range-optimal trajectories over 10000 ft
require higher airspeed, the second piece
should reach horizontal flight
at 10000 ft in order to be followed by an acceleration at 10000 ft.
\biglf
The climb at maximal admissible continuous thrust $T_{max}(h)$
and prescribed
airspeed $V_F(h)$ is completely determined by the initial conditions,
and \eref{eqTWfull} is solved\RSfootnote{singamma4fixedRTW.mws}
for $\gamma$ via
$$
2K\sin\gamma=R-\sqrt{R^2-4KRT/W+4KC_{D_0}R^2+4K^2}
$$
to get the flight path.
\biglf
Figure \RSref{Fig2p}\RSfootnote{twopiece.m} shows
such two-piece climbs at 250 KIAS, starting at 1500 ft
and stopping the first part at distance 1 to 4 nm in steps of 0.5.
The second part has an optiomally
reduced thrust and is stopped at $\gamma=0$.
For long-range flights, the trajectory ending at 10000 ft should be selected.
\begin{figure}[!]
\begin{center}
\includegraphics[width=4cm,height=4cm]{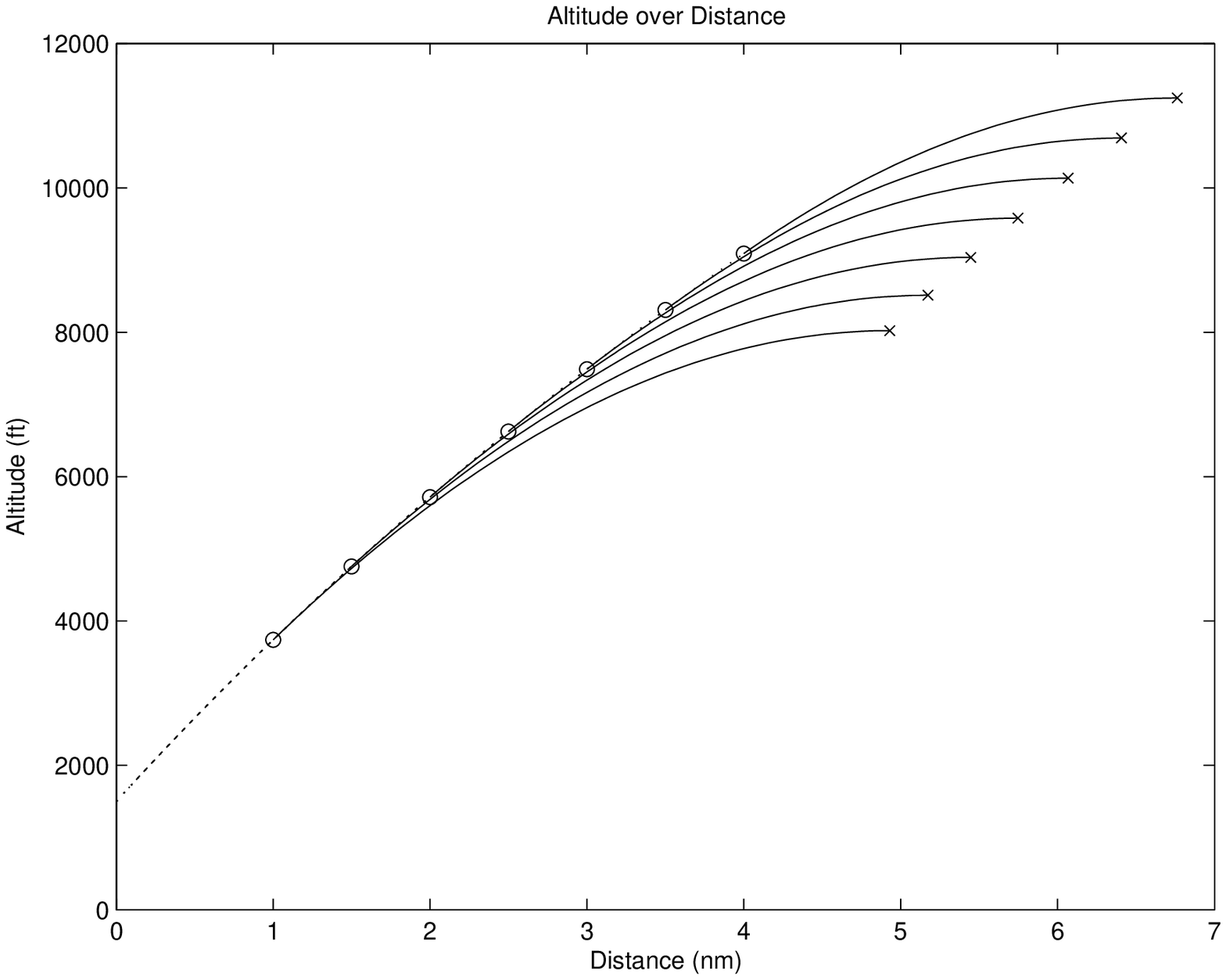}
\includegraphics[width=4cm,height=4cm]{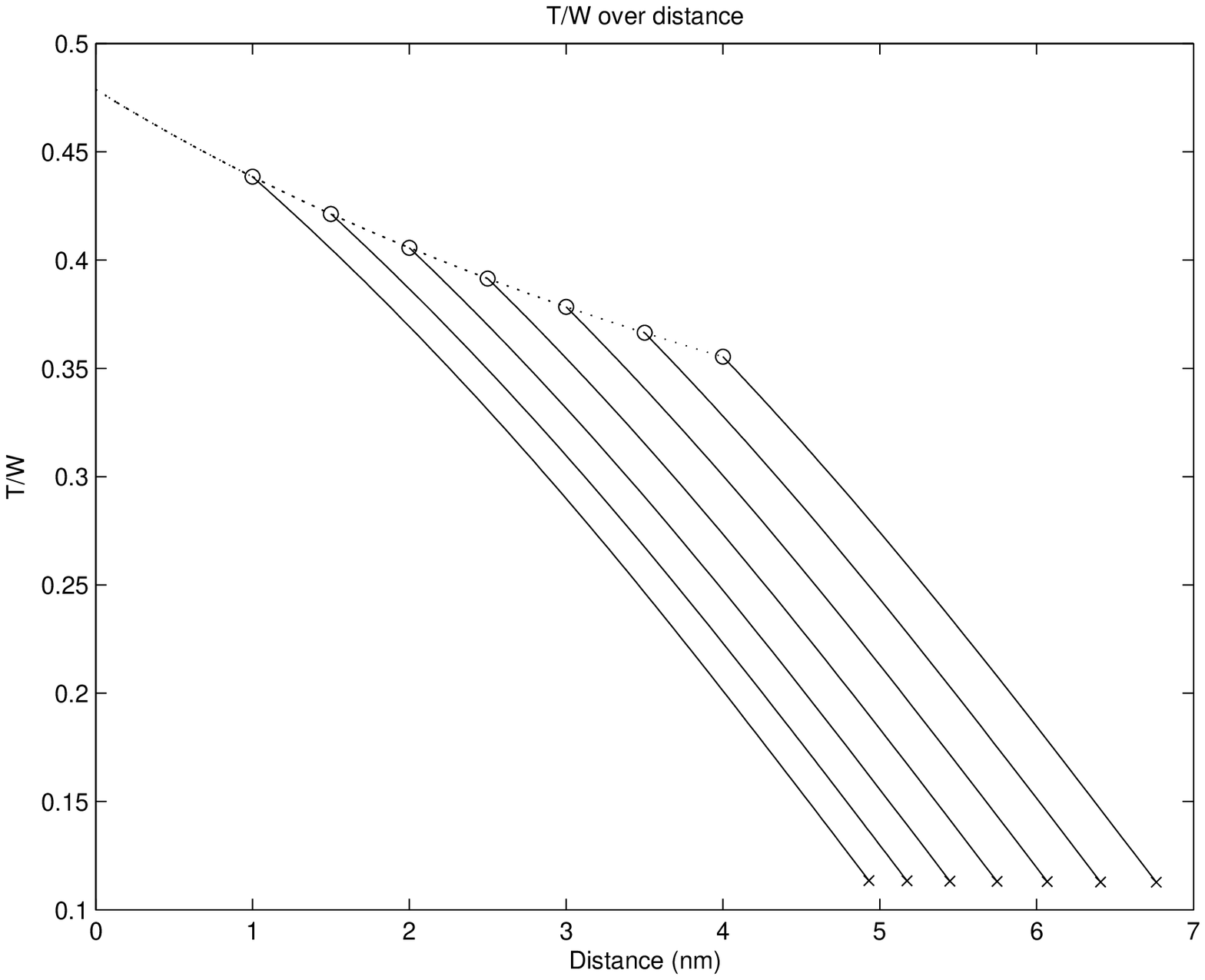}
\includegraphics[width=4cm,height=4cm]{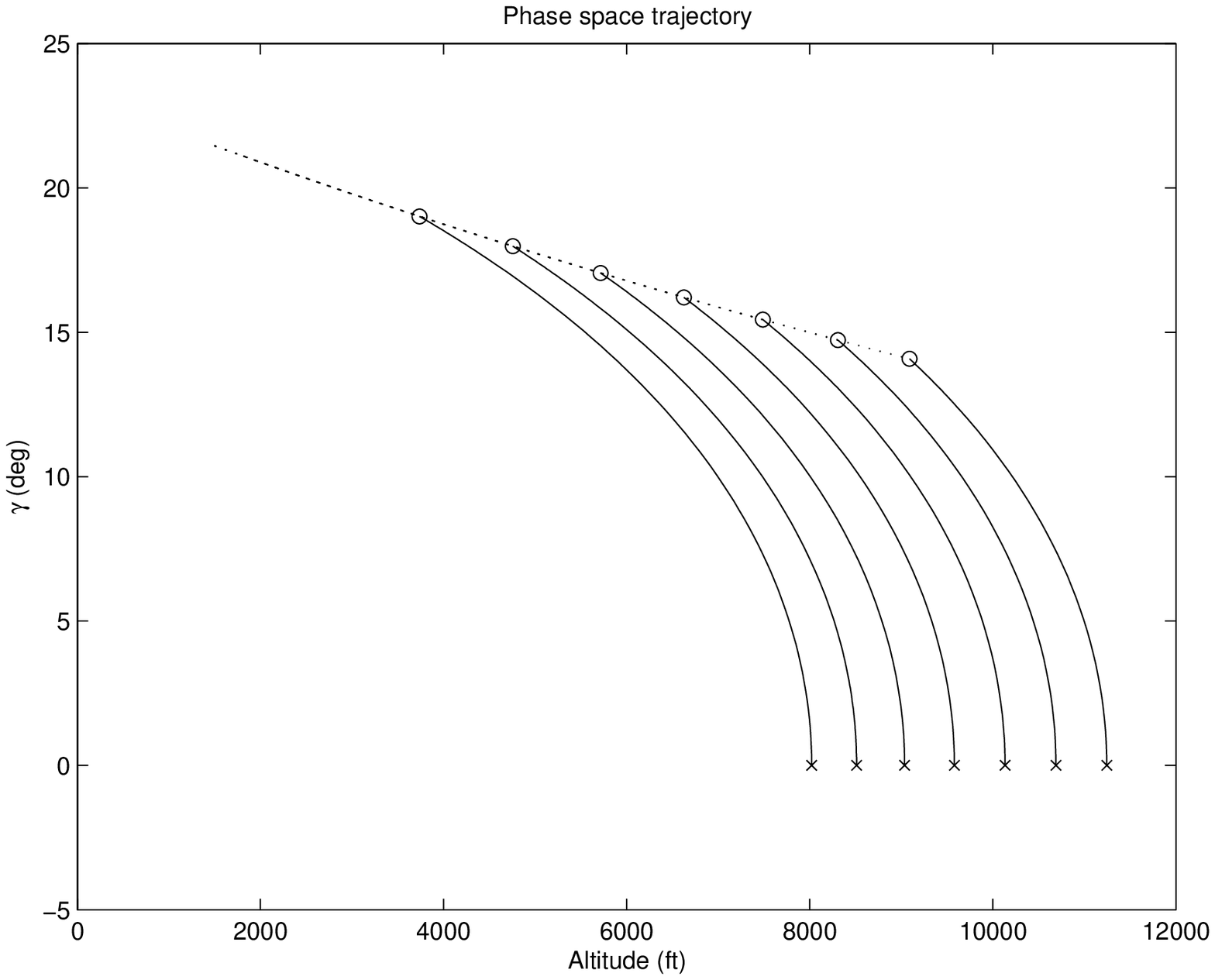}
\caption{Two-piece range-optimal solutions
at 250 KIAS,  flight path, $T/W$ ratio, and phase space 
   \RSlabel{Fig2p}}
\end{center} 
\end{figure}
\biglf
But these trajectories need 250 KIAS to be started,
and this calls for an
acceleration at the ``acceleration'' altitude where
clean configuration is reached and ``{\it at which the aircraft
accelerates towards the initial climb speed}'' \cite[p. 1245]{airbus:2011-1}.
Another acceleration will be necessary at 10000 ft,
because the range-optimal climb below 10000 ft is flown
at 250 KIAS, while the range-optimal climb to
higher altitudes starts at roughly 400 kts,
see Figure \RSref{Fig3pFP}, top right.
But if flown at high thrust, these two accelerations can be neglected
for long-range flights.
They take 1 nm and 6 nm,
respectively\RSfootnote{testAcceleration.m},
for the model aircraft of \RScite{hull:2007-1}.
%****************************************************************
\section{Flight Level Change in Cruise}\RSlabel{SecQSFFLC}
We now consider the practical situation that a long-distance
high-altitude cruise
under Air Traffic Control is a sequence of
level flights with various short-term flight-level changes.
These  are short-term changes
of $\gamma$, and
it is debatable whether they should be considered as quasi-steady flight.
We know now that such a flight is never
range-optimal, but each level section should
apply the $V_0$ speed given by \eref{eqV0}.
This means that all level flight
sections in cruise use the same $R_{0}$ from \eref{eqR0},
leading to the same 
$T/W$ ratio via \eref{eqTWfull},
no matter what the flight level or the propulsion model is.
Only the drag polar is relevant. Again, it turns out to
be convenient to work in terms of $R$ to be
independent of weight and altitude. 
\biglf
The $V_0$ speed at $\gamma=0$
then is a function of weight and altitude alone, and
flight level changes should comply with this, i.e. the speed should
still vary smoothly, while $\gamma$ and thrust may change rapidly.
We shall deal with this by keeping the
flight level change as quasi-steady flight,
except for the beginning and the end, where we allow an instantaneous
and simultaneous change of $\gamma$ and thrust that compensate each other.
\biglf
The idea is to keep the quasi-steady flight equation \eref{eqTWfull}
and the $R_{0}$ equation \eref{eqR0} valid at all times.
Then a jump in $\gamma$ must be counteracted by a jump in thrust,
one in the beginning and one in the end of the flight level change.
These instants are not quasi-steady, but the rest is.
\biglf
Consider a climb from altitude $h_0$ 
to altitude $h_1$. When flying
at $V_{0}$ at $h_0$ at maximal thrust $T_{max}(h_0)$,
the flight level change is impossible.
Otherwise, the quasi-steady flight equation \eref{eqTWfull}
at time $t_0$ and $\gamma=0$ is
\bql{eqT0W0}
\dfrac{T_0}{W_0}=C_{D_0}R_{0}+\dfrac{K}{R_{0}},\;T_0<T_{max}(h_0),
\eq
and we apply maximal thrust and go over to
$$
\dfrac{T_{max}(h_0)}{W_0}=C_{D_0}R_{0}+\dfrac{K}{R_{0}}\cos^2\gamma_0
+\sin\gamma_0
$$
defining a unique climb angle $\gamma_0$ satisfying
\bql{eqTransgam0}
2K\sin\gamma_0=R_{0}-\sqrt{R_{0}^2-4K
\dfrac{T_{max}(h_0)}{W_0}R_{0}+4KC_{D_0}R_{0}^2+4K^2}.
\eq
We could keep this angle for the climb, but we might reach the
thrust limit if we do so. Therefore we prefer to satisfy
$$
\dfrac{T_{max}(h)}{W(h)}=C_{D_0}R_{0}+\dfrac{K}{R_{0}}\cos^2\gamma(h)
+\sin\gamma(h)
$$
at each altitude using
$$
K\sin\gamma(h)=R_{0}-\sqrt{R_{0}^2-4K
\dfrac{T_{max}(h)}{W(h)}R_{0}+4KC_{D_0}R_{0}^2+4K^2}.
$$
This is put into an ODE system for $h$ and $W$ with $\gamma$
as an intermediate variable, namely
$$
\begin{array}{rcl}
  h'&=& \tan \gamma(h)\\
  W'&=& -\dfrac{C(h)T_{max}(h)}{V_{0}(h,W)\cos\gamma(h)}.
\end{array}
$$
The result is a climb with constant $R_{0}$ that keeps $V_{0}$
of \eref{eqV0} at all times and thus starts and ends with the
correct speed for range-optimal level flight. For descent, the same
procedure is used, but idle thrust is inserted. If the altitude change
is small, the solution is close to using the fixed climb/descent angle
$\gamma_0$ of
\eref{eqTransgam0}. At the end of the flight-level change at altitude
$h_1$,
the final speed $V_{0}(h_1,W_1)$ is the starting speed of the next level
flight, and the thrust has to be decreased instantaneously to
$T_1$ in order to keep the ratio
$$
\dfrac{T_0}{W_0}=\dfrac{T_1}{W_1}
$$
from \eref{eqT0W0}.
\biglf
We omit plots for our standard aircraft model, because they all
show that the crude simplification
$$
\dfrac{h_1-h_0}{x_1-x_0}\approx \gamma_0 \approx \dfrac{T-T_0}{W_0}
$$
holds for small altitude changes between level flights,
where the thrust $T$ is either $T_{max}$
or $T_{idle}$. Thus in $(x,h)$ space the transition
is very close to linear with the roughly constant climb angle
given above.
\biglf
But we have to ask whether climbing at maximal thrust is
fuel-to-distance optimal against all other choices of thrust.
If we insert the above approximation
into the fuel consumption with respect to the distance
and just keep the thrust varying, 
we get 
$$
  \int_{x_0}^{x_1}\dfrac{CT}{V\cos\gamma}dx
  \approx W_0(h_1-h_0)+T_0(x_1-x_0)
$$
up to a factor, and thus we should minimize the climb angle
if we relate consumption to distance. For descent, this leads to
taking $T_{idle}$ and is easy to obey, but for climb the range-optimal
solutions cannot be taken because they take too long.
Consequently, pilots are advised to
perform the climb at smallest rate allowed by ATC.
%****************************************************************
\section{Flight Phases for Maximal Range}\RSlabel{SecFPfOFU}
As long as Air Traffic Control does not interfere, we now see that a long
range-optimal flight should have the following phases:
\begin{enumerate}
\item Takeoff to clean configuration and acceleration altitude, 
\item accelerate there to 250 KIAS at maximal admissible continuous thrust,
%  using Section \RSref{SecAaFA},
\item climb at maximal admissible continuous thrust, keeping 250 KIAS and
  following the
  range-optimal angle selection strategy of Section \RSref{SecPS},
  and continuing with
  \item a solution of the variational problem given there
  to end at precisely 10000 ft in horizontal flight,
\item accelerate at 10000 ft in horizontal flight
  until the required speed
  for a range-optimal climb is reached,
%  using Section \RSref{SecAaFA} again,
\item perform a range-optimal climb/cruise following Sections
  \RSref{SecVP} and \RSref{SecRMaxHF} at maximal continuous admissible thrust
  until shortly before the top-of-descent point, leaving that climb for
\item an Euler-Lagrange path satisfying the variational problem
of Section \RSref{SecVP} until thrust is idle, 
\item do a continuous descent at idle thrust down to the Final Approach Fix.
\end{enumerate}
To arrive at the right distance and altitude, the time for
starting phase 7 needs to be be varied, like in Figure \RSref{Fig3pFP}.
\biglf
If ATC requires horizontal flight phases and  correspondent flight-level
changes, step 6 is followed
by
\begin{itemize}
\item[6a.] an Euler-Lagrange path satisfying the variational problem
  of Section \RSref{SecVP} to reach the prescribed altitude,
\item[6b.] using Section \RSref{SecRMaxHF} for
  range-optimal speed at level flight, and 
  \item[6c.] flight path changes following Section \RSref{SecQSFFLC},
\end{itemize} 
but the flight will not be
range-optimal.
Various examples show that a continuous descent from high altitude
ends up at speeds below 250 KIAS at 10000 ft, and deceleration is not needed.
\biglf
Flight paths for shorter distances should follow the above steps
for long-haul flights up to a certain
point
where they take a ``shortcut'' from the long-distance
flight pattern. 
%****************************************************************
\section{Conclusion}\label{SecOP}
Except for the two accelerations at 10000 ft and ``acceleration altitude'',
this paper provided range-optimal flight paths as simple solutions of certain
ordinary differential equations, without using Control Theory or other
sophisticated tools. However, everything was focused on
quasi-steady flight within simple atmosphere and propulsion models.
Also, the numerical examples were 
currently confined to the Small Business Jet of % Hull
\RScite{hull:2007-1}
with its turbojet engines.  
However, most of the results are general enough
to be easy to adapt for
other aircraft and engine characteristics, and this is left open.

%%%%%%%%%%%%%%%%%%%%%%%%
\bibliographystyle{plain}

\begin{thebibliography}{10}

\bibitem{airbus:2011-1}
Airbus.
\newblock {\em Airbus A 380 Flight Crew Operating Manual}.
\newblock Airbus S.A.S Customer Services Directorate, 31707 Blagnac, France,
  2011.
\newblock Reference: KAL A 380 Fleet FCOM, Issue Date: 03 Nov. 2011.

\bibitem{betts:1998-1}
J.T. Betts.
\newblock Survey of numerical methods for trajectory optimization.
\newblock {\em Journal of Guidance, Control, and Dynamics}, 21:193--207, 1998.

\bibitem{bryson-et-al:1969-1}
A.E. Bryson, M.N. Desai~Jr., and W.C. Hoffman.
\newblock Energy-state approximation in performance optimization of supersonic
  aircraft.
\newblock {\em Journal of Aircraft}, 6:481--488, 1969.

\bibitem{burrows:1982-1}
J.W. Burrows.
\newblock Fuel optimal trajectory computation.
\newblock {\em Journal of Aircraft}, 19:324--329, 1982.

\bibitem{calise:1977-1}
A.J. Calise.
\newblock Extended energy management methods for flight performance
  optimization.
\newblock {\em AIAA Journal}, 15:314--321, 1977.

\bibitem{franco-rivas:2014-1}
A.~Franco and D.~Rivas.
\newblock Analysis of optimal aircraft cruise with fixed arrival time including
  wind effects.
\newblock {\em Aerospace Science and Technology}, 32:212--222, 2014.

\bibitem{franco-et-al:2010-1}
A.~Franco, D.~Rivas, and A.~Valenzuela.
\newblock Minimum-fuel cruise at constant altitude with fixed arrival time.
\newblock {\em Journal of Guidance, Control, and Dynamics}, 33:280--285, 2010.

\bibitem{garciaheras-et-al:2016-1}
J.~García-Heras, M.~Soler, and F.J. Sáez.
\newblock Collocation methods to minimum-fuel trajectory problems with required
  time of arrival in {ATM}.
\newblock {\em Journal of Aerospace Information Systems}, 13:243--265, 2016.

\bibitem{gardi-et-al:2016-1}
A.~Gardi, R.~Sabatini, and S.~Ramasamy.
\newblock Multi-objective optimisation of aircraft flight trajectories in the
  {ATM} and avionics context.
\newblock {\em Progress in Aerospace Sciences}, 83:1--36, 2016.

\bibitem{huang-et-al:2012-1}
G.~Huang, Y.~Lu, and Y.~Nan.
\newblock A survey of numerical algorithms for trajectory optimization of
  flight vehicles.
\newblock {\em Sci. China Technol. Sci.}, 55:2538--2560, 2012.

\bibitem{hull:2007-1}
D.G. Hull.
\newblock {\em Fundamentals of Airplane Flight Mechanics}.
\newblock Springer, 2007.

\bibitem{maazoun:2015-1}
W.~Maazoun.
\newblock {\em Conception et analyse d'un système d'optimisation de plans de
  vol pour les avions}.
\newblock PhD thesis, École Polytechnique de Montréal, 2015.

\bibitem{miele:1962-1}
A.~Miele.
\newblock {\em Flight Mechanics: Theory of Flight Paths}.
\newblock Dover Books on Aeronautical Engineering, reprint of the 1962
  original, 2016.

\bibitem{myose-et-al:2005-1}
R.~Myose, T.~Young, and G.~Sim.
\newblock Comparison of business jet performance using different strategies for
  flight at constant altitude.
\newblock AIAA 5th ATIO and 16th Lighter-Than-Air Sys Tech. and Balloon Systems
  Conferences, 2005.

\bibitem{park-clarke:2015-1}
S.G. Park and J.-P. Clarke.
\newblock Optimal control based vertical trajectory determination for
  continuous descent arrival procedures.
\newblock {\em Journal of Aircraft}, 52:1469--1480, 2015.

\bibitem{peckham:1974-1}
D.H. Peckham.
\newblock {\em Range Performance in Cruising Flight}.
\newblock National Technical Information Service, 1974.

\bibitem{phillips:2010-1}
W.~F. Phillips.
\newblock {\em Mechanics of Flight}.
\newblock John Wiley \& Sons, 2010.
\newblock Second Edition.

\bibitem{pierson-ong:1989-1}
B.L. Pierson and S.Y. Ong.
\newblock Minimum-fuel aircraft transition trajectories.
\newblock {\em Mathematical and Computer Modelling}, 12:925--934, 1989.

\bibitem{rader-hull:1975-1}
J.E. Rader and D.G. Hull.
\newblock Computation of optimal aircraft trajectories using parameter
  optimization methods.
\newblock {\em Journal of Aircraft}, 12:864--866, 1975.

\bibitem{rutowski:1954-1}
E.S. Rutowski.
\newblock Energy approach to the general aircraft performance problem.
\newblock {\em Journal of the Aeronautical Sciences}, 21:187--195, 1954.

\bibitem{saucier-et-al:2017-1}
A.~Saucier, W.~Maazoun, and F.~Soumis.
\newblock Optimal speed-profile determination for aircraft trajectories.
\newblock {\em Aerospace Science and Technology}, 2017.
\newblock In Press.

\bibitem{stengel:2004-1}
R.F. Stengel.
\newblock {\em Flight Dynamics}.
\newblock Princeton University Press, 2004.

\bibitem{valenzuela-rivas:2014-1}
A.~Valenzuela and D.~Rivas.
\newblock Optimization of aircraft cruise procedures using discrete trajectory
  patterns.
\newblock {\em Journal of Aircraft}, 51:1632--1640, 2014.

\bibitem{vinh:1980-1}
N.X. Vinh.
\newblock {\em Optimal Trajectories in Atmospheric Flight}.
\newblock Elsevier, reprint in 2012, 1980.

\bibitem{vinh:1995-1}
N.X. Vinh.
\newblock {\em Flight Mechanics of High-Performance Aircraft}.
\newblock Cambridge University Press, 1995.

\bibitem{yajnik:1977-1}
K.S. Yajnik.
\newblock Energy-turn-rate characteristics and turn performance of an aircraft.
\newblock {\em Journal of Aircraft}, 14:428--433, 1977.

\end{thebibliography}

\end{document}